\documentclass{article}
\usepackage{epsfig}\usepackage{psfrag}
\usepackage{latexsym}
\usepackage{amsmath}
\usepackage{amsthm}
\usepackage{amssymb}
\usepackage{graphicx}
\usepackage{mdwlist}
\usepackage[all]{xypic}

\usepackage{hodge}

\begin{document}

\title{Gerby Localization, $\mathbb{Z}_3$-Hodge Integrals and the GW Theory of $\orb$}
\author{Charles Cadman\thanks{Supported by NSF grant No. 0502170.} and Renzo Cavalieri}

\maketitle
\begin{abstract}
We exhibit a set of recursive relations that completely determine all equivariant \GW\ invariants of $\orb$. 
We interpret such invariants as \zhh, and produce relations among them via Atiyah-Bott localization on moduli spaces of twisted stable maps to gerbes over $\proj$. 
\end{abstract}
\section*{Introduction}
\subsection*{Results}

Let $\mathbb{Z}_3\cong \mu_3=\{1,\omega,\bar\omega \}$ act on $\mathbb{C}^3$ via
$$
\omega \mapsto
\left[
\begin{array}{ccc}
 \omega & & \\
  & \omega & \\
  &  & \omega
\end{array}
\right]
$$
and denote $\orb$ the corrseponding stack quotient.
The natural linear action of the three dimensional torus $(\Cstar)^3$ descends to the quotient. We study the equivariant orbifold \GW\ invariants of this stack. \GW\ theory for orbifolds is developed in \cite{cr:ogwt} and \cite{cr:nctoo}. The algebraic point of view is established in \cite{Abramovich_Graber_Vistoli}.

\begin{mtheorem}
We exhibit a set of recursive relations that effectively\footnote{Effectively means that these recursions can be used to actually compute any invariant one may be interested in. The recursions have been coded to produce the invariants in section \ref{toi}.} compute any equivariant \GW\ invariant of $\orb$.  We also translate the recursions into a system of PDEs.
\end{mtheorem}

Following  \cite{bgp:crc}, we interpret such invariants as \zhh\ on moduli spaces of $\mathbb{Z}_3$-admissible covers, and deduce relations among them via Atyiah-Bott localization. Our computations, together with the associativity of the quantum cohomology of $\orb$, in fact allow us to evaluate a large family of such integrals.

Let $\adm{n_1}{n_2}$ be the component of the space of genus $0$ twisted stable maps to $\mathcal{B}\bZ_3$ corresponding to $\bZ_3$-admissible covers of an unparametrized $\bP^1$ with $n_1$ marked points having $\omega$-monodromy and $n_2$ having $\bar\omega$-monodromy (notation \ref{notation_A}).  There are no unramified marked points.

\begin{propp}\label{tpzhi}
Denote by $\lambda_i$ the $i$-th Chern class of the bundle $\mathbb{E}_\omega$ on $\adm{n_1}{n_2}$. All three-part \zhh
\begin{eqnarray}
	\int_{\adm{n_1}{n_2}} \lambda_i\lambda_j\lambda_k\nonumber
\end{eqnarray}
are effectively computed.
\end{propp}
\subsection*{History and Connections}
The orbifold $\orb$ has recently been an exciting object of study for both mathematicians and physicists.  In mirror symmetry, $\orb$ represents a distinguished point (the \textit{orbifold point}) in the $A$-model moduli space for the \textit{local $\Proj$} (the total space of the canonical bundle of $\Proj$). By studying topological strings on the mirror $B$-model, Mina Aganagic, Vincent Bouchard and Albrecht Klemm predicted several \GW\ invariants for $\orb$ (\cite{abk:tsaamf}).

The quest for verifying mathematically the physicists' predictions turned out to be more challenging than expected. After much effort, it was fulfilled almost simultaneously, and with independent methods, by three different ``teams''. Besides the current work, we acknowledge:
\begin{description}
\item[\cite{ccit:c3z3}:] Tom Coates, Alessio Corti, Hiroshi Iritani and Hsian-Hua Tseng apply Givental's formalism and an extension of the quantum Riemann-Roch and quantum Lefschetz theorems (\cite{c:qrrls}) to the orbifold setting. They compute the \textit{twisted J-function}, a generating function that encodes the invariants of $\orb$ after a change of variables (the inverse of the \textit{mirror map}). The inverse of the mirror map is not available in closed form, but its Taylor expansion can be computed to any order. This allows them to extract the invariants. Coates, Corti, Iritani and Tseng were the first to confirm the predictions of Aganagic, Bouchard and Klemm.
\item[\cite{bc:c3z3}:] Arend Bayer and the first author found an explicit way to  construct the space of $n$-pointed, genus $0$ twisted stable maps to $B\mu_r$ from $\overline{M}_{0,n}$ using $r$-th root contructions.  Then they used the formalism of weighted stable maps to write down a new and explicit formula for the total Chern class of the obstruction bundle for Gromov-Witten invariants of $[\bC^N/\mu_r]$.  From this, they derived a combinatorial formula for the non-equivariant \GW\ invariants of $\orb$.
\end{description}

The possibility of understanding and developing connections between these three different approaches is by itself an exciting perspective. Our interest in the \GW\ theory of $\orb$ is further enhanced by the possibility of framing its study into a broader context. In particular, we briefly discuss the connections with a fascinating conjecture of Yongbin Ruan (CRC), and with the study of tautological classes on the moduli space of curves.

\subsubsection*{The Crepant Resolution Conjecture}

Mirror symmetry and the general philosophy of the McKay correspondence motivated Yongbin Ruan to formulate the following conjecture.

\vspace{0.2cm}
\noindent \textbf{\textit{Ruan's Crepant Resolution Conjecture} \cite{r:crc}.} \\ \textit{The quantum cohomology rings of a Gorenstein\footnote{A Gorenstein orbifold is an orbifold having generically trivial stabilizers whose canonical bundle is pulled back from a line bundle on its coarse moduli space.  The terminology is standard and should not be confused with ``Gorenstein stack.''} orbifold $\mathfrak{X}$ and of any crepant resolution $Y$ (if one exists) are isomorphic.} 
\vspace{0.2cm}

In 2005 (\cite{bg:crc}), Jim Bryan and Tom Graber verified, Ruan's conjecture in some examples and proposed a very strong reformulation of the conjecture: the \GW\ potentials of $\mathfrak{X}$ and $Y$ should be equal after a linear transformation on the cohomology insertion variables and the specialization to roots of unity of the excess quantum parameters (see \cite[Section 1.3]{bg:crc} for the precise statement). Several more examples are checked in \cite{bgp:crc}, \cite{betal}, \cite{davesh}, \cite{ccit:crcas}.

Coates, Corti, Iritani and Tseng remark that in all examples verified a technical condition on the orbifold cohomology holds (\textit{hard Lefschetz}, \cite[Definition 1.1]{bg:crc}). They supply evidence that  Bryan-Graber's formulation of the CRC should be modified when hard Lefschetz does not hold, and propose the conjecture should be phrased in terms of Givental's formalism. All  \GW\ invariants for a space $X$ are encoded in the geometry of  a Lagrangian cone $\mathcal{L}_X$ inside Givental's  symplectic vector space $\mathcal{H}_X=H^\ast(X) \otimes \bC((z^{-1}))$.

\vspace{0.2cm}
\noindent \textbf{\textit{Crepant Resolution Conjecture} \cite{ccit:crcs}.} \\ \textit{There is a degree preserving  $\bC((z^{-1}))$- linear symplectic isomorphism $\mathbb{U}\colon \mathcal{H}_{\fX}\rightarrow \mathcal{H}_{Y}$  such that, after analytic continuation, we have $\mathbb{U}(\mathcal{L}_{\fX})=\mathcal{L}_Y$.} \vspace{0.2cm}

Further, $\mathbb{U}$ satisfies three assumptions that we choose not to record here. This point of view is very powerful and may be the correct avenue to get to a general proof of the CRC.  However, we feel it worthwhile to seek a direct description of the relationship between the \GW\ invariants of $\fX$ and $Y$. 
In other words, an attractive question is: \textit{what is the strongest and most simple-minded formulation of the CRC that holds in general?}
The orbifold $\orb$ is the first meaningful example where the hard Lefschetz condition does not hold. We hope understanding its \GW\ theory will help to answer this question.

\subsubsection*{\zhh}

The Hodge bundle $\mathbb{E}$ is a rank $g$ vector bundle on $\overline{\mathcal{M}}_g$ whose fibers over any curve are the global sections of the dualizing sheaf. In \cite{m:taegotmsoc}, David Mumford proves that the Chern classes $\lambda_i=c_i(\mathbb{E})$ are tautological, and suggests their study as an approach to intersection theory on the moduli space of curves. Using \GW\ theory as a key tool, Carel Faber and Rahul Pandharipande (\cite{f:algo}, \cite{fp:lsahiittr}, \cite{fp:hiagwt}) carry on Mumford's program, and unveil beautiful structure underlying Hodge integrals: appropriate generating functions for Hodge integrals are governed by the classical KdV hierarchy.
Further connections with Hurwitz theory are established by the celebrated ELSV formula (\cite{elsv:ohnahi}, \cite{elsv:hnaiomsoc}) and are used in interesting work 
of Ravi Vakil with Tom Graber (\cite{gv:taut}, \cite{gv:rvl}) and  with Ian Goulden and David Jackson (\cite{gjv:gwpophnahh}, \cite{gjv:spolgc}, \cite{gjv:last}) making progress towards a combinatorial description of the tautological ring of the moduli space of curves.

We  turn our attention to moduli spaces of curves that admit a finite group action. In this case the Hodge bundle splits into eigenbundles corresponding to the decomposition of a fiber into irreducible representations. The Chern classes of such subbundles give rise to a new set of tautological classes worth investigating. Such classes are extremely well behaved from an intersection theoretic point of view: 
they ``split naturally'' along the boundary, and satisfy natural ``Mumford-type'' relations analogous to those for ordinary Hodge classes. Following a suggestion of Ruan, we call top intersections of such classes $G$\textit{-Hodge integrals}. First, a natural curiosity is whether $G$-Hodge integrals are naturally encoded in some natural integrable system. Second, by studying $G$-Hodge integrals we intend to strengthen the connection between the geometry of the moduli space of curves and representation/Hurwitz theory, in the hope of gaining insight towards the structure of the tautological intersection ring for the moduli space of curves.

\subsection*{Strategy and Techniques}
Our approach to the study of $G$-Hodge integrals is similar in spirit to Faber and Pandharipande's \cite{fp:hiagwt}.
They recognize Hodge integrals in the localization fixed loci contributions of some auxiliary integral on the very degenerate moduli space ${\overline{M}_{g}(\proj,1)}$. This produces a wealth of relations among Hodge integrals that can subsequently be inverted.

In the setting of $G$-Hodge integrals one  looks at twisted stable maps to $\mathcal{B}G\times\bP^1$, which can be viewed as admissible $G$-covers of $\bP^1$, and applies the localization formula for a 1-dimensional torus action on the base $\bP^1$.  This is exploited by the second author, Aaron Bertram and Gueorgui Todorov (\cite{bct:gg-1}) to give a purely combinatorial proof of a celebrated Faber-Pandharipande formula (the computation of $\lambda_g\lambda_{g-1}$ on the hyperelliptic locus). 

For $\orb$, localization on spaces of admissible covers does not seem to provide enough relations to compute the \zhh; therefore we looked for a new idea.  If one thinks of these admissible covers as stable maps to $\mathcal{B}\bZ_3\times\bP^1$, it is natural to replace $\mathcal{B}\bZ_3\times\bP^1$ with a similar looking stack.  Such a stack $\fG$ should have the following two properties.
\begin{enumerate}
\item There is a morphism $\fG\to\bP^1$ whose fibers are isomorphic to $\mathcal{B}\bZ_3$.  This ensures that stable maps to $\fG$ locally look like admissible covers.
\item The inertia stack of $\fG$ is isomorphic to $\fG\times\bZ_3$.  This ensures that the monodromy around a fixed point of the $\bZ_3$-action is well-defined (see section~\ref{glocal}).
\end{enumerate}
The stacks which have these properties are called $\bZ_3$-gerbes over $\bP^1$.  There are two such stacks up to isomorphism\footnote{There are three $\bZ_3$-gerbes over $\bP^1$, but the nontrivial ones are isomorphic as stacks.}:  $\mathcal{B}\bZ_3\times\bP^1$ and a nontrival gerbe which we denote $\fG_1$.  By localizing on a space of twisted stable maps to $\fG_1$, we found the relations we needed to compute all the \zhh.

In slightly more detail, our strategy for computing  \zhh\  combines two types of relations:
\begin{description}
	\item[WDVV:] the associativity of quantum cohomology provides a set of relations that allows to express any integral on a space with ``many" $\bar\omega$ points in terms of integrals on spaces with at most $2$ $\bar\omega$ points. Such integrals are to be considered as initial conditions.
	\item[Localization:] the evaluation via localization of auxiliary integrals on spaces of  maps to gerbes  provide a set of recusions among \zhh\ on spaces with at most  $2\  \bar\omega$ points. This determines all initial conditions in terms of the three pointed integrals, which can be computed by hand.
	\end{description}
	
 Auxiliary integrals on  moduli spaces of stable maps to a $\mathbb{Z}_3$-gerbe over $\proj$ must have the following characteristics.
\begin{enumerate}
	\item The integral often vanishes for dimension reasons.  When it doesn't, we can use two different linearizations of the vector bundles involved to get a nontrivial relation (the integral is independent of the choice of linearization).
	\item When the integral is evaluated via localization, the contributions of the various fixed loci contain \zhh. 
	\item The combinatorial complexity of the fixed loci contributions can be kept under control. In section \ref{restr} we explain a few ``tricks'' used to achieve this goal.
\end{enumerate}
This strategy allows one to produce a massive amount of relations between \zhh. Somewhat surprisingly, a large number of relations are (non-trivially) dependent. Only after much effort did we obtain enough relations to completely determine the \GW\ theory of $\orb$.

\begin{rmk} In sections \ref{locrel} and \ref{proof} we present our relations in a form that makes it easiest to prove how they inductively compute all the invariants of $\orb$. The data of infinitely many recursive relations is efficiently (if more obscurely) packaged in a handful of differential equations on  appropriate generating functions (section \ref{gf}).
\end{rmk}

\subsection*{Plan of the Paper}

Section \ref{glocal} is devoted to the developement of the technique of gerby localization. We assume a little familiarity with the \GW\ theory of stacks as in \cite{Abramovich_Graber_Vistoli}. We  give a fairly extensive working presentation of gerbes, and of Atyiah-Bott localization in the context of maps to gerbes. 

In section \ref{gworb} we discuss the \GW\ theory of $\orb$. In particular, we describe all equivariant invariants in terms of \zhh and show that WDVV imposes strong conditions on these invariants.

Section \ref{locrel}  carefully develops the localization computations that  produce  relations between \zhh. It is hard to avoid being technical with such computations. We sought transparency by adding comments and building the contributions  via  ``elementary" pieces  that are subsequently organized in tables. 

After the hard work of section \ref{locrel}, proving the Main Result in section \ref{proof} is  a short and leisurely stroll.

Section \ref{gf} is a ``commercial'' for the language of generating functions. Here we translate the information of infinitely many messy recursions into a handful of PDE's between appropriate generating functions.

Finally, in section \ref{toi} we  include a significant number of invariants of $\orb$, in case some sharp eye could help us detect some structure. It would be extremely nice to have a closed form description of the potential of $\orb$. 
\subsection*{Acknowledgements}
We would like to warmly acknowledge Jim Bryan for suggesting this problem to us; Rahul Pandharipande for boosting our morale by telling us we were attempting a ``very difficult computation''. 
Tom Coates for providing us with a nice exposition of their statement of the crepant resolution conjecture. Dan Abramovich, Aaron Bertram, Max Lieblich, Melissa Liu, Sam Payne and Michael Thaddeus for helpful conversations. Aarend Bayer for his priceless proofreading and for noticing the relationship between our generating functions and the $J$ function. Last but not least Y.P. Lee who sat with us through a very long afternoon during which the details of gerby localization started to take concrete shape.

 \section{Gerby localization}
\label{glocal}
\subsection{Background and motivation}

Let $\fX$ be a smooth Deligne-Mumford stack having a projective coarse moduli scheme $X$. In \cite{Abramovich_Vistoli}, Abramovich and Vistoli defined twisted stable maps to $\fX$ and showed that the connected components form proper Deligne-Mumford stacks. They come with natural evaluation maps to (a rigidification of) the inertia stack of $\fX$, which we denote $\sI(\fX)$.  There is a natural perfect obstruction theory on the stack of twisted stable maps, and therefore the standard algebraic definition of Gromov-Witten invariants \cite{Behrend} works if one uses insertions coming from the cohomology of $\sI(\fX)$.

A twisted stable map to $\fX$ over a scheme $S$ is a commutative diagram
\begin{equation}
\label{tw_st_map}
\xymatrix{\Sigma_i \ar@{^{(}->}[r] \ar[d] & \fC \ar[r] \ar[d] \ar@/_/[dd] & \fX \ar[d] \\
\sigma_i \ar@{^{(}->}[r] \ar[dr]_{\cong} & C \ar[d] \ar[r] & X \\
 & S,}
\end{equation}
where:
\begin{itemize}
	\item $C\to X$ is an ordinary stable map over $S$, with sections $\sigma_i$;
	\item $\fC$ is a twisted curve with coarse moduli space $C$;
	\item $\fC\to\fX$ is representable;
	\item $\Sigma_i$ (the markings) are \'{e}tale gerbes over $\sigma_i$. 
\end{itemize}
	 These gerbes can be constructed by applying root constructions along $\sigma_i\subseteq C$.  See \cite{Cadman_root} for more on root constructions.  The morphism $\fC\to C$ is an isomorphism away from the gerbes $\Sigma_i$ and the singular locus of $\fC\to S$.  The fibers of $\fC\to S$ can have twisted nodal singularities.

For a finite group $G$, letting $\fX=\mathcal{B}G$ leads to a theory of $G$-covers of curves, which was studied in \cite{Abramovich_Corti_Vistoli}.  In this case one doesn't need as much stack machinery, since the moduli problem can be defined in terms of certain $G$-covers of nodal curves.  Given an action of $G$ on $Y$, one can similarly compare twisted stable maps to a quotient stack $[Y/G]$ with equivariant stable maps to $Y$, which was used in \cite{Jarvis_Kaufmann_Kimura}.

\subsection{Gerbes}

We recall the definition of a gerbe from \cite[3.15]{Laumon_Moret-Bailly}.

\begin{Def}
A \textbf{gerbe} over a scheme $X$ is a stack $\fX$ equipped with an epimorphism $\fX\to X$ such that the diagonal $\fX\to\fX\times_X \fX$ is an epimorphism.
\end{Def}

In other words, a stack $\fX$ over $X$ is a gerbe if local sections exist and if any two local sections are locally isomorphic.  It is perplexing at first to think of sections as being isomorphic, but the categorical nature of stacks is essential for this definition.  If any two sections were locally ``equal'', then they could be glued to give a global section.  However, not all gerbes have global sections, and the existence of gerbes without global sections is essential for the calculations done in this paper.  One can think of a gerbe as a sheaf of categories (cf. \cite[Def. 2.1]{Moerdijk}).

If $X=Spec\; k$, with $k$ an algebraically closed field, then any gerbe over $X$ is isomorphic to $\mathcal{B}G$ for some finite group $G$.  If $X$ is a $k$-variety, then there is a trivial gerbe with fiber $\mathcal{B}G$ over $X$, namely $X\times \mathcal{B}G$.  The fibers of a gerbe can be ``twisted'' in at least two different ways.  One is by an element $\xi\in H^1(X,Aut(G))$.  Such an element determines a fiber bundle with fiber $G$ and structure group $Aut(G)$, which is the same as a group scheme over $X$.  Associated to this group scheme is its classifying stack, which one might naively regard as the fiber bundle with fiber $\mathcal{B}G$ and structure group $Aut(G)$ associated to $\xi$.

Suppose now that $G$ is an abelian group; let $\tilde{G}\to X$ be the group scheme associated to $\xi$ as above, and let $B\tilde{G}$ be its classifying stack.  The inertia stack of $\mathcal{B}\tilde{G}$ is then isomorphic to $\tilde{G}\times_X \mathcal{B}\tilde{G}$.  This shows that the nontrivial variation of the group $G$ over $X$ can be detected by the inertia stack of $\mathcal{B}\tilde{G}$, even though the stack $\mathcal{B}\tilde{G}$ naively seems not to contain that information (the $k$-points of $\mathcal{B}\tilde{G}$ are the same as the $k$-points of $X$).  This brings us to the definition of a $G$-gerbe, which encapsulates the other way in which the fibers of a gerbe can be ``twisted.''

\begin{Rem}
For nonabelian groups $G$, the definition of $G$-gerbe is more complicated than the one below.  The notions of gerbe and band were formulated by Giraud \cite{Giraud}.
\end{Rem}

\begin{Def}
Let $G$ be a finite abelian group.  A gerbe with \textbf{band} $G$ is a gerbe $\fX\to X$ together with an isomorphism of group stacks over $\fX$, $\sI(\fX)\to G\times\fX.$  We demand that the following diagram 1-commute.
$$\xymatrix{\sI(\fX) \ar[r] \ar[rd] & G\times\fX \ar[d] \\
 & \fX.}$$
\end{Def}

This is the same as defining, for each object $x$ of $\fX$, an isomorphism $G\to Aut(x)$ which is compatible with restrictions and isomorphisms of objects of $\fX$.  Gerbes with band $G$ over $X$ are classified by $H^2(X,G)$.

If $n$ is a positive integer, then from the exact sequence $$1\to \mu_n \to \bG_m \to \bG_m \to 1,$$ there is a homomorphism $H^1(X,\bG_m)\to H^2(X,\mu_n)$.  In other words, to any line bundle $L$ on $X$, one can associate a $\mu_n$-gerbe over $X$.  The total space of this $\mu_n$-gerbe is the stack $X_{(L,n)}$ defined as follows.

\begin{Def}
A section of $X_{(L,n)}\to X$ over a morphism of schemes $f:S\to X$ is a pair $(M,\varphi)$, where 
\begin{enumerate}
\item $M$ is a line bundle on $S$ and 
\item $\varphi:M^{\otimes n}\to f^*L$ is an isomorphism.
\end{enumerate}
A morphism in the category $X_{(L,n)}$ from $(M,\varphi)$ to $(N,\psi)$ over a commutative diagram
$$\xymatrix{S \ar[r]^h \ar[dr]_f & T \ar[d]^g \\ & X}$$
is an isomorphism $\rho:M\to h^*N$ such that the following diagram commutes.
$$\xymatrix{M^{\otimes n} \ar[r]^{\rho^n} \ar[d]_{\varphi} & h^*N^{\otimes n} \ar[d]^{\psi} \\
f^*L \ar[r]^{\sim}_{\text{canonical}} & h^*g^*L}$$
\end{Def}

Here are some important facts about $X_{(L,n)}$.
\begin{enumerate}
\item This has the structure of $\mu_n$-gerbe if $\mu_n$ acts on each object $(M,\varphi)$ of $X_{(L,n)}$ by multiplication on $M$.
\item An isomorphism $\sigma:L_1\otimes P^{\otimes n}\to L_2$ of line bundles on $X$ induces an isomorphism of $\mu_n$-gerbes $X_{(L_1,n)}\to X_{(L_2,n)}$ sending $(M,\varphi)$ to $(M\otimes P,\sigma\circ(\varphi\otimes 1))$.
\item On $\bP^1$, the $\mu_n$-gerbes are classified by $\bZ_n$, with $k$ corresponding to the gerbe $\bP^1_{(\sO_{\proj}(k),n)}$.
\item If $P$ is the complement of the zero section in $L$, then $X_{(L,n)}$ is isomorphic to $[P/\bC^*]$, where $\bC^*$ acts as the $n$-th power of the standard action \cite[2.3.5]{Cadman_root}.
\end{enumerate}

\subsection{Atiyah-Bott Localization Formula}

We give a brief account of localization and develop some details geared to our application of it. In our  treatment we follow the ``localization language'' and notations in \cite[chapters $4$ and $27$]{clay:ms}.

Consider the  one-dimensional algebraic torus $\mathbb{C}^\ast$, and recall that the $\mathbb{C}^\ast$-equivariant Chow ring of a point is a polynomial ring in one variable:
$$A^\ast_{\Cstar}(\{pt\},\mathbb{C})= \mathbb{C}[\hbar]. $$

Let $\Cstar$ act on a smooth, proper Deligne-Mumford stack $X$, denote by $i_k:F_k\hookrightarrow X$ the irreducible components of the fixed locus for this action and by $N_{F_k}$ their normal bundles. The natural map:
$$
\begin{array}{ccc}
A^\ast_{\Cstar}(X) \otimes_{\mathbb{C}[\hbar]} \mathbb{C}(\hbar) & \rightarrow & \sum_{k}{A^\ast_{\Cstar}}(F_k) \otimes_{\mathbb{C}[\hbar]} \mathbb{C}(\hbar)\\
                                             &             &                                        \\
\alpha                                       & \mapsto     &\displaystyle{\sum_k\frac{i_k^\ast\alpha}{c_{top}(N_{F_k})}}.
\end{array}
$$
is an isomorphism. Pushing forward equivariantly to the class of a point, we obtain the Atiyah-Bott integration formula:
$$\int_{[X]}\alpha = \sum_k \int_{[F_k]} \frac{i_k^\ast\alpha}{c_{top}(N_{F_k})}.$$
The extension of this formula to smooth Deligne-Mumford stacks was established by \cite{Graber-Pandharipande}.

To illustrate a confusing aspect of localization on stacks, consider the following example.  Let $\fX$ be the square root of $\bP^1$ at $0$ and let $\bC^*$ act on $\bP^1$ in the standard way fixing $0$ and $\infty$.  By choosing a linearization of $\sO_{\proj}(1)$ having weight $0$ at $\infty$, one can lift this action to $\fX$ (such liftings are explained in the following subsection).  We'll write $\iota$ for inclusions of fixed loci.  If $\sigma$ is the Poincare dual of $0$ in $A^1_{\bC^*}(\bP^1)$, then localization on $\bP^1$ gives $$\int_{\bP^1}\sigma = \int_{(\bP^1)^{\bC^*}}\frac{\iota^*\sigma}{\hbar} = 1,$$  since $\iota^*\sigma$ is $\hbar$ at $0$ and $0$ at $\infty$.  If one does this on $\fX$, then the fixed locus over $0$ is a copy of $\mathcal{B}\bZ_2$.  We identify the Chow group of $\mathcal{B}\bZ_2$ with that of a point by pullback, so that integration introduces a factor of $1/2$.  Then localization on $\fX$ gives
$$\int_{\fX} \sigma = \int_{\mathcal{B}\bZ_2} \frac{\hbar}{\hbar/2} = 1.$$  In order to get the right answer, we had to put in a normal bundle weight of $\hbar/2$.  The fractional weight is essentially due to the fact that the tangent sheaf of $\fX$ is not a $\bC^*$-linear sheaf in the sense of \cite[4.3]{Romagny}.  One must first pass to a $\bZ_2$-extension of $\bC^*$.  However, $\bC^*$ has an honest action on the tangent bundle to $\fX$, as a stack, induced by the one on $\fX$.  Alternatively, it is because the restriction of the $\bC^*$-action to $\mathcal{B}\bZ_2$ is not the same as the trivial action, though each element of $\bC^*$ acts trivially.  As remarked in \cite[5.3]{Kresch}, one needs the $\bZ_2$-extension of $\bC^*$ in order to trivialize the action on $\mathcal{B}\bZ_2$.

There is a generalization of the above example to $r$-th root constructions.  It also affects the computations done in this paper when looking at the normal direction to a fixed locus in a space of twisted stable maps which smooths a twisted node.  This introduces a factor of $1/r$ relative to the smoothing of the node on the coarse moduli space.

\subsection{Maps to $\bZ_3$-gerbes over $\bP^1$}

  Let $\fG_i=\bP^1_{(\sO_{\proj}(i),3)}$.  We will only consider the cases $i=0,1$.  
  
\begin{Rem} 
  $\fG_0=\bP^1\times B\bZ_3=[\bP^1/\bZ_3]$, with $\bZ_3$ acting trivially. Twisted stable maps to $\fG_0$ are admissible covers of a parametrized $\bP^1$ (\cite{r:dd2}). 
\end{Rem}

 Given a twisted stable map to $\fG_k$, the twisted marked points can be separated into $\omega$ points and $\bar{\omega}$ points.  To make this distinction, we first need a canonical generator of the automorphism group of a twisted marking on a twisted curve.  The twisted curve is locally the quotient of a cyclic cover which is totally ramified at the marked point.  An oriented simple loop around the branch point determines an element of the Galois group of the cover, which is the canonical generator of the stabilizer group of the fixed point.

Recall that $\mu_3$ acts compatibly on all objects of $\fG_i$, and therefore each point of $\fG_i$ has its automorphism group identified with $\mu_3$.  Given a twisted stable map $f:\fC\to \fG_i$ and a twisted point $x\in\fC$, $x$ is called an $\omega$-point if $f_*:Aut(x)\to\mu_3$ sends the canonical generator to $\omega$.  Otherwise it is an $\bar{\omega}$-point. Note that the generator cannot be sent to $1$, since $f$ is representable.
\begin{Rem}
 This definition of $\omega$ and $\bar{\omega}$-points can be rephrased in terms of evaluation maps.  The evaluation map at an $\omega$-point maps to the $\omega$-component of the (rigidified) inertia stack, and likewise for $\bar{\omega}$-points.
\end{Rem}
 
We denote by
\begin{eqnarray*}
	\sG_i(k,\ell) \subseteq \overline{M}_{0, k+\ell}(\fG_i,1)
\end{eqnarray*}
 the component of the space of $k+\ell$-marked, degree $1$, genus $0$, twisted stable maps to $\fG_i$ parameterizing maps that have $k$ $\omega$-points and $\ell$ $\bar{\omega}$-points.  We measure degree by composing with the map $\fG_i\to\bP^1$. 
There is a universal diagram
$$
\begin{array}{ccccc}
\mathcal{G}_i(k,\ell)_1 & \stackrel{f}{\longrightarrow} & \fG_i  & \longrightarrow & \proj \\
 & & & &  \\
 \pi \downarrow & & & & \\ 
 & & & &  \\
\mathcal{G}_i(k, \ell), & & & &
\end{array}
$$
where $\mathcal{G}_i(k,\ell)_1$ is the universal curve over $\mathcal{G}_i(k,\ell)$.

\subsection{Localization Set-up}
In order to apply localization to spaces of twisted stable maps to gerbes, we first show that they are smooth.

\begin{Thm} \label{thm:smooth_dimension}
Let $\fX$ be an \'etale gerbe over a homogeneous space $X$ (i.e., $X=G/P$ for some semi-simple complex Lie group $G$ and parabolic subgroup $P$, and $\fX\to X$ is \'etale).  Then each connected component of the stack of genus $0$ twisted stable maps to $\fX$ is smooth.  Moreover, the natural map $\overline{\sM}_{0,n}(\fX) \to \overline{\sM}_{0,n}(X)$ sends a connected component to a component of the coarse moduli space having the same dimension.
\end{Thm}

\begin{Prf}
For smoothness, it suffices to show that for any genus $0$ twisted stable map $F:\fC\to\fX$, $H^1(\fC,F^*T_{\fX})=0$.  Let $f:C\to X$ be the associated map of coarse moduli spaces, and denote the natural maps $\pi:\fX\to X$ and $\nu:\fC\to C$.  Since $f^*T_X$ has nonnegative degree on each component of $C$, it follows that $H^1(C,f^*T_X)=0$.  Moreover, $\pi^*T_X=T_{\fX}$ since $\fX\to X$ is \'etale.  Therefore, $F^*T_{\fX}=\nu^*f^*T_X$.  Since $\nu$ has $0$-dimensional fibers and $\nu_*\sO_{\fX}=\sO_X$, it follows that $$H^i(\fC,F^*T_{\fX})=H^i(C,f^*T_X)$$ for $i=0,1$.  This implies smoothness and shows that the tangent spaces of $\overline{\sM}_{0,n}(\fX)$ and $\overline{\sM}_{0,n}(X)$ have the same dimension at the points determined by $F$ and $f$.
\end{Prf}

\vspace{0.5cm}
Recall from above that $\fG_i$ is a global quotient $[P/\bC^*]$, where $P$ is the complement of the zero section in the total space of $\sO_{\proj}(i)$, and $\bC^*$ acts as the cube of the standard action.  Therefore, to define a $\bC^*$ action on $\fG_i$, it suffices to define a $\bC^*$ action on $\bP^1$ together with a linearization of $\sO_{\proj}(i)$.  We choose the linear action of $\bC^*$ on $\bP^1$ fixing $0,\infty$ and acting with weight $1$ on the tangent space at $0$ and $-1$ at $\infty$.  For the linearization, we choose weights $i$ at $0$ and $0$ at $\infty$.

We digress slightly to formulate our group action in terms of Definition 2.1 of \cite{Romagny}, which helps clarify the need for the linearization of $\sO_{\proj}(i)$.  We need to define a morphism $\bG_m\times \fG_i\to \fG_i$.  We denote the action of $\bG_m$ on $\bP^1$ defined above by $(t,f)\mapsto tf$.  Suppose we are also given a linearization of $\sO_{\proj}(i)$, so that for any $f:S\to \bP^1$ and $t:S\to\bG_m$ we have an isomorphism $\alpha_t:f^*\sO_{\proj}(i)\to (tf)^*\sO_{\proj}(i)$, such that $\alpha_u\circ\alpha_t = \alpha_{ut}$.  Then we can define $\bG_m\times \fG_i\to \fG_i$ by $$(t,f,M,\varphi)\mapsto (tf,M,\alpha_t\circ\varphi).$$  It is easy to verify that this defines a \emph{strict} action, a notion defined in [ibid].

The action on $\fG_i$ induces an action on $\sG_i(k,\ell)$ by post composition.  The equivariant vector bundles on $\sG_i(k,\ell)$ we use in our localization integrals come from line bundles on $\fG_i$ by pulling back to the universal curve over $\sG_i(k,\ell)$ and taking the $R^1$-pushforward. 

There are two ``types'' of line bundles on $\fG_i$, that generate $Pic(\fG_i)$:
\begin{itemize}
	\item line bundles pulled back from $\proj$, which we denote $\mathcal{O}(n)$;
	\item a tautological cube root of $\sO_{\proj}(i)$.
\end{itemize}
  For $i=0$, we denote this cube root by $L_{\omega}$, because it is identical to a trivial bundle on $\bP^1$ on which $\bZ_3$ acts by $\omega=e^{2\pi i/3}$.  On $\fG_1$, the tautological cube root is denoted $\sO(1/3)$.  At each point of $\fG_1$, $\sO(1/3)$ also restricts to the $\omega$ representation of $\bZ_3$.  
 \begin{Rem} 
  While $L_{\omega}$ can be linearized trivially, $\sO(1/3)$ has degree $1/3$, and therefore the weights of the linearization at $\infty$ and $0$ must differ by $1/3$.
\end{Rem}
Let us consider the pullback of the tautological cube root to a twisted curve $\fC$, where $f:\fC\to \fG_i$ is a degree $1$ twisted stable map having $k$ $\omega$-points and $\ell$ $\bar{\omega}$-points.  As the map has degree $1$, the pullback has degree $i/3$.  We can compute the fractional part of the degree by noting that an $\omega$-point contributes $1/3$ to the fractional part and an $\bar{\omega}$-point contributes $2/3$.  We therefore arrive at the relation
\begin{eqnarray}
	k-\ell\equiv i \text{\ mod\ } 3.
	\label{gerbymonodormy}
\end{eqnarray}
The spaces $\sG_i(k,\ell)$ are empty when the above relation does not hold.

\vspace{0.5cm}
Now we compute the Euler characteristic of the pullback of a line bundle from $\fG_i$.  We use the Riemann-Roch formula for twisted curves from \cite[7.2.1]{Abramovich_Graber_Vistoli}:
$$\chi(\sE) = rk(\sE)\chi(\sO_{\fC}) + \mbox{deg}(\sE) - \sum_{j=1}^{k+\ell} \mbox{age}_{p_j}(\sE),$$ where $\sE$ is a vector bundle on a twisted cuves $\fC$ and $p_1,\ldots,p_{k+\ell}$ are the twisted points. 

 On $\fG_0$, let $\sL=\sO(n)\otimes L_{\omega}$, and let $f:\fC\to \fG_0$ be a map in $\sG_0(k,\ell)$.  The age of $f^*\sL$ at an $\omega$-point is $1/3$, and at an $\bar{\omega}$-point it is $2/3$.  Therefore, we have
\begin{equation}\label{Euler_char_1}
\chi(f^*(\sO(n)\otimes L_{\omega}))=n+1-\frac{k+2\ell}{3}.
\end{equation}  

On $\fG_1$, for $\sL=\sO(n+1/3)$:
\begin{equation}\label{Euler_char_2}
\chi(f^*\sO(n+\frac{1}{3})) = n + \frac{4}{3} - \frac{k+2\ell}{3}.
\end{equation}

Two types of localization integrals appear in this paper.  The first set have the form
\begin{equation}\label{localization_integral_1}
\int_{\sG_i(k,\ell)}e(R^1\pi_*f^*\sE)\cup \rho
\end{equation}
where $e$ denotes the equivariant Euler class, $\sE$ is a rank $3$ equivariant vector bundle on $\fG_i$, and $\rho$ is a product of classes of the form $ev_j^*(0)$ and $ev_j^*(\infty)$.  Since $0$ and $\infty$ are torus fixed points and $ev_j$ is equivariant, $\rho$ is an equivariant cohomology class.  The vector bundle $\sE$ is a direct sum of three line bundles, and we choose various weights for their linearizations at $0$ and $\infty$, subject to the restriction that the weight at $0$ minus the weight at $\infty$ equals the degree of the line bundle.

\begin{Prop}
Assume $k+\ell>0$ and let $\sE$ be one of the following bundles.
\begin{enumerate}
\item $(\sO\oplus\sO\oplus\sO(-1))\otimes L_{\omega}$ on $\fG_0$
\item $\sO(-2/3)\oplus\sO(-2/3)\oplus\sO(-2/3)$ on $\fG_1$
\end{enumerate}
Then the integral in (\ref{localization_integral_1}) is $0$ whenever $$\ell-c_1(\sE)+\deg(\rho)< 3.$$
\end{Prop}

\begin{Prf}
We claim that the pullback of $\sE$ by any map in $\sG_i(k,\ell)$ has vanishing global sections.  In the case $i=1$, this is clear since each summand has negative degree on every component.  For $i=0$, one only needs to check that the pullback of $L_{\omega}$ has no non-trivial sections.  Since it has degree $0$, the only way it could have a section is if it were trivial.  But this is impossible if either $k$ or $\ell$ is positive, since every section would have to vanish at a twisted point.

From Theorem \ref{thm:smooth_dimension}, we have $$\dim \sG_i(k,\ell)=k+\ell,$$ so the integral vanishes whenever $rk(R^1\pi_*f^*\sE)+\deg(\rho)<k+\ell$.  From the formulas (\ref{Euler_char_1},\ref{Euler_char_2}), we see that $$rk(R^1\pi_*f^*\sE)=k+2\ell-c_1(\sE)-3.$$
\end{Prf}

The other set of integrals we consider have the form
$$\int_{\sG_0(p,q)}\lambda_{top}\lambda_{top-i}\lambda_{top-j}\cup ev_1^*(0)\cup ev_2^*(0)\cup ev_3^*(\infty),$$
with $(p,q)=(3k+3,0), (3k+1,1)$ or $(3k-1,2)$. 

  Here $\lambda_n=c_n((R^1\pi_*f^*L_{\omega})^{\vee})$, where $L_{\omega}$ is given the trivial linearization (with weights $0$). On the  moduli spaces considered, $R^1\pi_*f^*L_{\omega}$ is a bundle of rank $k$.  We also use $\lambda_i$ for the analogous class on the space of admissible covers.

\subsection{ Fixed Loci Contributions}
\label{restr}

\begin{figure}[hbt]
	\begin{center}
	\psfrag{0}{$0$}
	\psfrag{i}{$\infty$}
	\psfrag{l}{$(k_0,\ell_0)$}
	\psfrag{r}{$(k_\infty,\ell_\infty)$}
	\psfrag{C}{$\fC$}
	\psfrag{1}{$1$}
	\psfrag{G}{$\fG_i$}
		\includegraphics[width=0.98\textwidth]{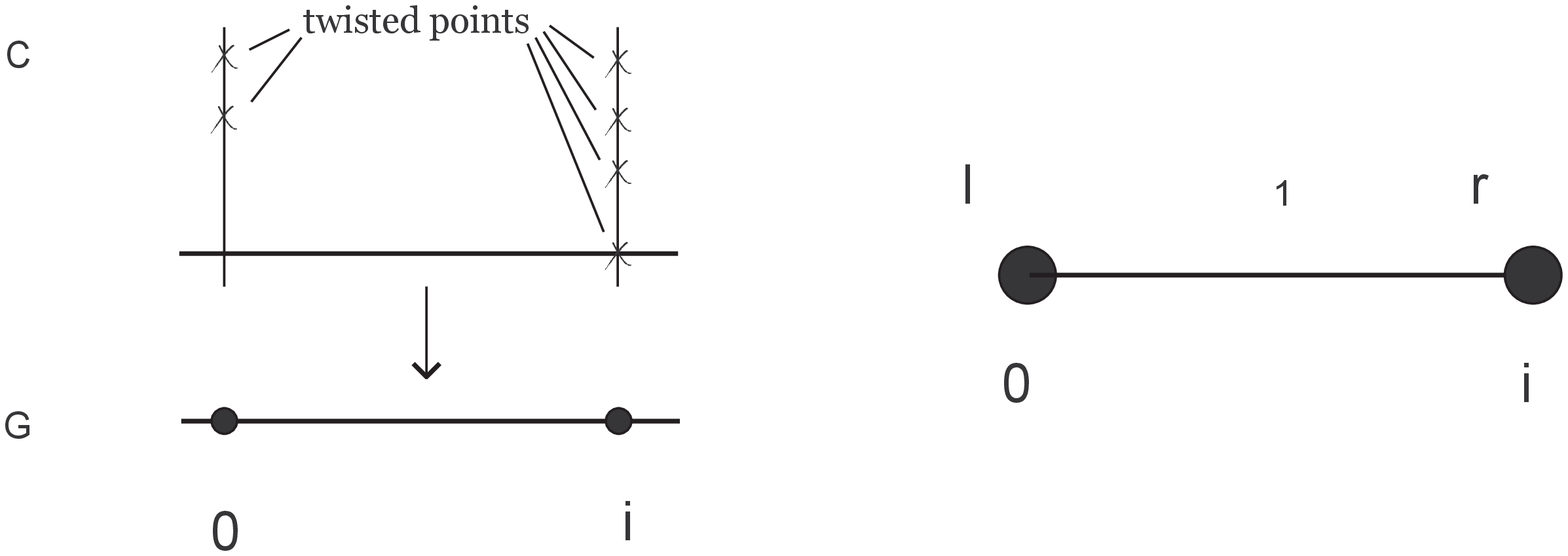}
	\end{center}
	\caption{A fixed locus for maps of degree $1$ to a $\proj$-gerbe and its associated localization graph.}
	\label{fix}
\end{figure}

The fixed loci for the induced action on the moduli space consist of maps such that anything ``interesting'' (branching, collapsed components, twisted points and marked points) happens over $0$ or $\infty$. Restricting our attention to maps of degree $1$, we have a main component mapping with degree $1$ to the gerbe, and possibly two contracted components over $0$ and $\infty$. The nodes can be twisted. All other marks and twisted points are on the contracted twigs. This is illustrated in Figure \ref{fix}, where we also show the associated localization graph.
	
For any such localization graph, there are several corresponding fixed loci given by all possible labellings of the marked points on the contracted components; this contributes a combinatorial factor.  The restriction of the equivariant cohomology class to the fixed locus can be analyzed by applying the normalization sequence to the bundle $f^*\sE$:
$$0\to f^*\sE\to \nu_*\nu^*f^*\sE \to \bigoplus_{\text{nodes\ } n} f^*\sE_n \to 0.$$
Here $\nu$ is the normalization map.  Since we will choose $f^*\sE$ to have no sections on any component, we will have an associated short exact sequence:
$$0\to \oplus H^0(f^*\sE_n) \to H^1(f^*\sE) \to H^1(\nu_*\nu^*f^*\sE) \to 0.$$  So
$R^1\pi_*f^*\sE$ splits into a sum in $K$-theory:
$$R^1\pi_*f^*\sE = \sum_n R^0\pi_*f^*\sE_n + \sum_i R^1\pi_*(f^*\sE)|_{C_i},$$
where the last sum is over the components of the curve $\fC$. We denote the main component by $C_0$  and by $C_1, C_2$ the contracted components. 

For the main component $C_0$, the Euler class of $R^1\pi_*(f^*\sE)|_{C_0}$ is called the \textbf{edge contribution}.  This bundle is trivial over the fixed locus, so it only contributes a weight factor.

To compute this weight factor we first choose the open cover $C_0=U_0\cup U_{\infty}$, where $U_0$ is the complement of $0$ and $U_{\infty}$ is the complement of $\infty$.  Let $\sL$ be a line bundle on $C_0$ and let $x$ be the coordinate on $\bP^1$ (the coarse moduli space of $C_0$) at $0$.  Then we can write
$$H^1(C_0,f^*\sE) = \Gamma(U_0\cap U_{\infty},f^*\sE)/(s_0\bC[x]+s_{\infty}\bC[x^{-1}]),$$ where $s_0$ is a \emph{minimally vanishing section} of $\sL$ on $U_0$ and $s_{\infty}$ is such a section on $U_{\infty}$.  This means that $s_0$ is nonvanishing away from $0$ and vanishes at $0$ to the lowest order possible, which is determined by the age of $\sL$ at $0$.  We consider two examples.
\smallskip

\noindent {\bf Example 1.}  Suppose that $i=0$, that $C_0$ has an $\omega$-point at $0$ and an $\bar{\omega}$-point at $\infty$, and that $\sL = f^*\sO(-1)\otimes L_{\omega}$.  This is a degree $-1$ bundle on $C_0$, and due to the twisting, $s_0$ vanishes to order $1/3$ at $0$ while $s_{\infty}$ vanishes to order $2/3$.  Therefore, $s_{\infty} = x^{-2}s_0$ and $H^1$ is generated by $x^{-1}s_0$, which has a pole of order $2/3$ at $0$.  Since the weight of $\bC^*$ on the tangent space at $0$ is $1/3$, the weight of this section is $\alpha+2/3$, where $\alpha$ is the weight of $\sL$ at $0$.  If we had interchanged $\omega$ and $\bar{\omega}$, then the weight would be $\alpha+1/3$.
\smallskip

\noindent {\bf Example 2.}  Suppose that $i=1$, that $C_0$ has $\bar{\omega}$-points at $0$ and $\infty$, and that $\sL=f^*\sO(-2/3)$.  Then both $s_0$ and $s_{\infty}$ vanish to order $2/3$ and $s_{\infty}=x^{-2}x_0$.  We again have that $H^1$ is generated by $x^{-1}s_0$, and the weight is $\alpha+1/3$, with $\alpha$ being the weight of $\sL$ at $0$.
\smallskip

The fiber of $f^*\sE$ at a node can only contribute sections if it has an eigenspace on which the stabilizer group of the node acts trivially.  Due to our choices for $\sE$, this will only happen when the node is untwisted, and in this case $R^0\pi_*f^*\sE_n$ will contribute a product of three weight factors determined by the linearizations. 

\begin{Rem}
If at least one of the summands of $\sE$ has weight $0$ at $0$ and $\infty$, then it follows from the preceding that any fixed locus having an untwisted node on the main component $C_0$ contributes $0$ to the localization formula.
\end{Rem}
The Euler class of the normal bundle to a fixed locus has a contribution from each vertex. There is a pure weight factor (``moving of the point'') and a term of the form (weight - $\psi$) (``smoothing of the node''). This part is standard; we refer to \cite[Chapter 27]{clay:ms} for details.  The \textbf{vertex contribution} at $C_i$ is the quotient of the Euler class of $R^1\pi_*(f^*\sE|C_i)$ by the normal bundle contribution from this vertex.

One subtle issue in the normal bundle contribution is that deforming a twisted node to first order does not deform the node on the coarse curve.  In other words the normal spaces are different, so the torus weights are not exactly the same.  In this paper, we always take our $\psi$-classes to be those living on the coarse curve, so when we smooth a twisted node, the factor will be (weight-$\psi/3$).  There is also a factor of $1/3$ in ``weight'' relative to the coarse moduli space.  For example, if we are deforming a twisted node at $0$, where the torus has weight $\hbar$ on the tangent space, then the factor will be $(\hbar/3-\psi/3)$.

Since this factor is placed in the denominator, the upshot is a factor of $3$ in the numerator for each twisted node.  There are also factors of $1/3$ which appear when comparing the fixed locus in the moduli space to the ``abstract fixed locus,'' where we forget $C_0$ and look at the contracted components, viewed as spaces of admissible covers.  To derive these factors, we need the following lemma.

\begin{Lem}
Let $\fX$ and $\fY$ be Deligne-Mumford stacks, let $H$ be an abelian group, and let $\fX\to \mathcal{B}H$ and $\fY\to\mathcal{B}H$ be morphisms.  Then the canonical morphism $\fX\times_{\mathcal{B}H}\fY\to\fX\times\fY$ is surjective, finite, and \'etale of degree $|H|$.
\end{Lem}

\begin{Prf}
It is a general fact that we have a fiber square
$$\xymatrix{\fX\times_{\mathcal{B}H} \fY \ar[r] \ar[d]\ar@{}[dr]|{\Box} & \fX\times\fY \ar[d] \\
\mathcal{B}H \ar[r] & \mathcal{B}H\times\mathcal{B}H.}$$
So the lemma follows from the fact that for abelian groups, the multiplication map $H\times H\to H$ is a group homomorphism, and so we have a fiber square
$$\xymatrix{\mathcal{B}H\ar[r] \ar[d]\ar@{}[dr]|{\Box} & \mathcal{B}H\times\mathcal{B}H \ar[d]\\
Spec\;\bC \ar[r] & \mathcal{B}H.}$$
\end{Prf}

We can apply this to a more general situation where $G$ is a finite group and $\fG$ is a $G$-gerbe over $\bP^1$.  Let $\sM_1$ and $\sM_2$ be two spaces of twisted stable maps into $\fG$ which send a particular marked point to $0$.  Suppose that the monodromy action around this marked point is $g\in G$ at one of the points and $g^{-1}$ at the other.  Let $H$ be the quotient of the centralizer of $g$ in $G$ by the subgroup generated by $g$.  Then the marked points can be glued, and the stack parametrizing the glued maps is isomorphic to $\sM_1\times_{\mathcal{B}H}\sM_2$ 
\cite[5.2]{Abramovich_Graber_Vistoli}.  The previous lemma shows that this differs from $\sM_1\times\sM_2$ by a factor of $|H|$.

If $\sM_1\times_{\mathcal{B}H}\sM_2$ were a fixed locus, then for the purpose of localizing, we view this as $|H|$ copies of $\sM_1\times\sM_2$.  However, as we mentioned prior to the lemma, there is an extra factor equal to the order of $g$ which enters into localization via the smoothing of the twisted node.  Therefore, the node overall contributes a factor equal to the order of the centralizer of $g$.

This analysis easily extends to multiple nodes, with each node contributing such a factor.  In our situation, we want to forget about the main component $C_0$ and focus solely on the contracted components.  Therefore, we have to divide by $3$ for the $\mu_3$ automorphisms on $C_0$.  In the end, our fixed loci pick up a \textbf{gluing factor} of $3^{N-1}$, where $N$ is the number of nodes.

In summary, for each localization graph as in Figure \ref{fix}, we have to consider
\begin{enumerate}
\item a combinatorial factor coming from the many ways to distribute marked points,
\item a vertex contribution which combines the restriction of the vector bundle with the Euler class of the normal bundle,
\item an edge contribution which comes from $H^1$ on the main component,
\item weight factors at untwisted nodes which sometimes kill off fixed loci,
\item and a gluing factor of $3^{N-1}$, where there are $N$ nodes.
\end{enumerate}



\section{\GW\ Theory of $\orb$}
\label{gworb}

Orbifold \GW\ invariants of $\orb$ are intersection numbers on the moduli space of twisted stable maps $\sm{n}{\orb}{\beta}$, of classes pulled-back from $H^\ast_{orb}(\orb)$ (i.e. the cohomology of the inertia stack).

\vspace{0.2cm}
\noindent\textbf{Observations:}
\begin{enumerate}
	\item The orbifold $\orb$ contains no compact curve classes, therefore the only invariants correspond to constant maps ($\beta=0$).
	\item The inertia stack $\mathcal{I}\orb$ is a disconnected union of $\orb$ and two copies of $\mathcal{B}\mathbb{Z}_3$, whose fundamental classes we identify with third roots of unity. Keeping track of the age grading, we obtain the three-dimensional ring:
	$$
	H^\ast_{orb}(\orb)=H^0\oplus H^2\oplus H^4= \mathbb{C}_{1} \oplus \mathbb{C}_{\omega} \oplus \mathbb{C}_{\bar\omega}.
	$$ 
\end{enumerate}

We denote the general (primitive) \GW\ invariant:
\begin{multline}\label{inv}
	\langle 1^{n_0} \omega^{n_1} \bar\omega^{n_2} \rangle = \\ \int_{[\sm{n_0+n_1+n_2}{\orb}{0}]^{vir}}\prod_1^{n_0}ev_i^\ast(1)\prod_1^{n_1}ev_j^\ast(\omega)\prod_1^{n_2}ev_k^\ast(\bar\omega)
\end{multline}
\subsection{Orbifold Invariants and $\mathbb{Z}_3$-Hodge Integrals} \label{oizhh}

For $n_1+n_2 \geq 3$, consider a \GW\ invariant $\langle\omega^{n_1} \bar\omega^{n_2} \rangle$.\footnote{For the sake of lighter notation, we omit from this discussion invariants with fundamental class insertions. Such invariants only appear in the three-pointed case and are discussed in section \ref{small}}

This invariant is supported on components of $\sm{n_1+n_2}{\orb}{0}$ parameterizing maps from curves that are non-trivially twisted at the marked points. All such maps must factor through the image of $\underline{0}\in \mathbb{C}^3$, the unique stacky point in $\orb$. Therefore, by introducing the euler class of an obstruction bundle to compare the (virtual) fundamental classes on the two different moduli spaces, the invariant can be computed as an integral over the moduli space $\sm{n_1+n_2}{\bz}{0}$:

\begin{eqnarray}
	\langle \omega^{n_1}\bar\omega^{n_2} \rangle =\int_{[\sm{n_1+n_2}{\bz}{0}]}e(Ob)\prod_1^{n_1}ev_i^\ast(\omega)\prod_1^{n_2}ev_j^\ast(\bar\omega)
	\label{bzinv}
\end{eqnarray}
 
In \cite{Abramovich_Corti_Vistoli}, Abramovich, Corti and Vistoli show that the stack $\sm{n_1+n_2}{\bz}{0}$ is the (normalization of the) moduli space of admissible $\mathbb{Z}_3$-covers of a genus $0$ curve. Informally, this stack parameterizes degree $3$ covers $p:E\rightarrow C$ such that:
\begin{itemize}
	\item $C$ is a stable $(n_1+n_2)$-marked genus zero curve;
	\item $E$ is a nodal curve and nodes of $E$ ``correspond to''\footnote{The preimages of nodal (resp. smooth) points of $C$ are nodal (resp. smooth) points of $E$.} nodes of $C$;
	\item $E$ is endowed with a $\mathbb{Z}_3$ action;
	\item $p$ is the quotient map with respect to the action;
	\item $p$ is ramified only over the marked points of $C$, and possibly over the nodes;
	\item when $p$ is ramified over a node, denote $x_1$ and $x_2$ the shadows of the node in the normalization $\tilde{E}$. The $\mathbb{Z}_3$-representations induced on $T_{x_1}$ and $T_{x_2}$ are dual to each other. 
\end{itemize}

The cohomology class $ev_i^\ast(\omega)$ (resp. $ev_i^\ast(\bar\omega)$) corresponds to selecting components (of the moduli space parameterizing covers) where the local monodromy around the $i$-th mark coincides with (resp. is dual to) the representation on the tangent space at the preimage of the mark.

\begin{notation} \label{notation_A} We denote by 
$$\mathcal{A}(n_1,n_2)$$
the component of $\sm{n_1+n_2}{\bz}{0}$ identified by the class 
$$  \prod_1^{n_1}ev_i^\ast(\omega)\prod_1^{n_2}ev_j^\ast(\bar\omega).$$
If, in addition, we have $m$ untwisted moving marked points, we adopt the notation
$$\mathcal{A}(n_1,n_2)_m.$$
\end{notation}
\begin{Rem}
 Since the total monodromy of a ramified cover of $\proj$ is $1$, we see that $\mathcal{A}(n_1,n_2)$ is non-empty only when
$$
n_1+2n_2 \equiv 0 \ \ (\mbox{mod}\  3).
$$
\end{Rem}

When the monodromy condition is verified, $\mathcal{A}(n_1,n_2)$ is a smooth stack of dimension $n_1+n_2-3$, with coarse moduli space $\overline{M}_{0, n_1+n_2}$. 

By  Riemann-Hurwitz, the genus of the covers parameterized in $\mathcal{A}(n_1,n_2)$ is $g=n_1+n_2-2$. The natural forgetful morphism
$$
\mathcal{A}(n_1,n_2)\longrightarrow \overline{\mathcal{M}}_g
$$
allows to define a Hodge bundle $\mathbb{E}$ on admissible covers by pull-back (for more details see \cite{r:dd2}, section 1.3).

The obstruction bundle (see \cite[Section 3]{bgp:crc})
$$Ob= R^1\pi_\ast f^\ast(L_\omega \oplus L_\omega \oplus L_\omega)$$

can be described in terms 
of the Hodge bundle $\mathbb{E}$ on  $\mathcal{A}(n_1,n_2)$. The Galois action on the covers induces a $\mathbb{Z}_3$ action on $\mathbb{E}$, which gives a decomposition 
$$
\mathbb{E}= \mathbb{E}_1\oplus \mathbb{E}_\omega \oplus \mathbb{E}_{\bar\omega}
$$
 into eigenbundles (with respect to the action of the primitive generator of the group).
 
From \cite{bgp:crc}, section $3$\footnote{They use the formulation $(\mathbb{E}^\vee)_{\bar{\omega}}$, which is equivalent: the $\bar\omega$ eigenbundle of the dual of the Hodge bundle is in fact the dual to the $\omega$ eigenbundle of the Hodge bundle.}:
$$
R^1\pi_\ast f^\ast(L_\omega)= (\mathbb{E}_\omega)^\vee 
$$

The bundle $\mathbb{E}_\omega$ has rank $r=\frac{n_1+2n_2}{3}-1$.

\begin{notation} Since in this work we only use the Chern classes for $\mathbb{E}_\omega$, to avoid useless proliferation of subscripts, we denote: 
\begin{eqnarray*}
	c_i(\mathbb{E}_\omega) := \lambda_i.
\end{eqnarray*}
\end{notation}

\GW\ invariants are now expressed as $\mathbb{Z}_3$-Hodge integrals:

$$
	\begin{imp}
	\begin{array}{lcl}
	\displaystyle{\langle \omega^{n_1} \bar\omega^{n_2} \rangle } & 
	\displaystyle{=} &
	\displaystyle{\int_{[\mathcal{A}(n_1,n_2)]} e(\mathbb{E}^\vee_\omega\oplus\mathbb{E}^\vee_\omega\oplus\mathbb{E}^\vee_\omega)=}
	\\ & & \\
	& 
	\displaystyle{=} &
\displaystyle{	(-1)^{n_1+n_2-3}\sum_{i+j+k=n_1+n_2-3}t_1^{r-i}t_2^{r-j}t_3^{r-k}\int_{[\mathcal{A}(n_1,n_2)]} \lambda_i\lambda_j\lambda_k.
	}
	\end{array}
	\end{imp}
	$$
\begin{Rem}
Note that Proposition \ref{tpzhi} in the introduction follows immediately from this formula together with our Main Result.
\end{Rem}

\subsection{Three-pointed Invariants}\label{small}
In the previous section we did not discuss invariants with fundamental class insertions. It is an  easy consequence of the projection formula that such invariants vanish if they contain more than three insertions. The three-pointed invariants are as follows.  Note that integration over $[\bC^3/\bZ_3]$ is defined by the localization formula.
\begin{description}
	\item[$\langle1^3\rangle$=] $\displaystyle{\int_{[\sm{3}{\orb}{0}]^{vir}}ev_1^\ast(1)\cup ev_2^\ast(1)\cup ev_3^\ast(1)=\int_{\orb} 1 =\frac{1}{3t_1t_2t_3}}$. \\
		\item [$\langle1 \omega\bar\omega\rangle$=]$\displaystyle{\int_{[\sm{3}{\orb}{0}]^{vir}}ev_1^\ast(1)\cup ev_2^\ast(\omega)\cup ev_3^\ast(\bar\omega)=\int_{\mathcal{A}(1,1)_1}1=\frac{1}{3}}$. 
		\item [$\langle \omega^3\rangle$=]$\displaystyle{\int_{[\sm{3}{\orb}{0}]^{vir}}ev_1^\ast(\omega)\cup ev_2^\ast(\omega)\cup ev_3^\ast(\omega)=\int_{\mathcal{A}(3,0)}1=\frac{1}{3}}$. 
		\item [$\langle \bar\omega^3\rangle$=]$\displaystyle{\int_{[\sm{3}{\orb}{0}]^{vir}}ev_1^\ast(\bar\omega)\cup ev_2^\ast(\bar\omega)\cup ev_3^\ast(\bar\omega)=\int_{\mathcal{A}(0,3)}\prod_1^3(\lambda_1+t_i)=\frac{t_1t_2t_3}{3}}$. 
\end{description}

\subsection{WDVV Relations}
\label{WDVVrel}
WDVV equations encode the associativity of the quantum product. One can think of them as infinitely many relations between the \GW\ invariants of $\orb$, or as a unique PDE on the \GW\ potential of $\orb$. In this section we develop the first point of view, which leads to an immediate proof of Proposition \ref{wdvvrel}. 
Consider the map
$$
\varphi: \mathcal{A}(n_1+2,n_2+2)\rightarrow \overline{M}_{0,4}=\proj,
$$
that records only the location of two $\omega$ and two $\bar\omega$ branch points.
The integral
\begin{eqnarray}
	\int_{[\mathcal{A}(n_1+2,n_2+2)]} e(Ob)\cup \varphi^\ast([pt])
	\label{wdvveq}
\end{eqnarray}
is independent of the choice of a representative for the class of a point in $\proj$. Equating the explicit evaluations for
$$
\begin{array}{ccccccccc}
P_1 =  & 
	\parbox{1cm}{	\includegraphics[width=0.15\textwidth]{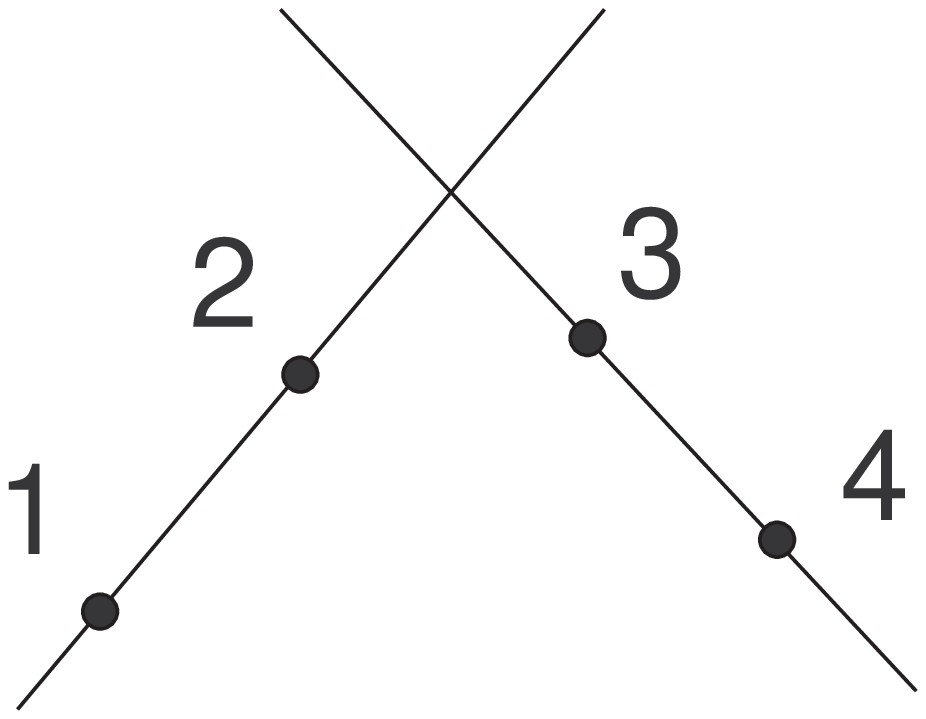}}
& \in \overline{M}_{0,4} & & ,& &P_2 =  & 
	\parbox{1cm}{	\includegraphics[width=0.15\textwidth]{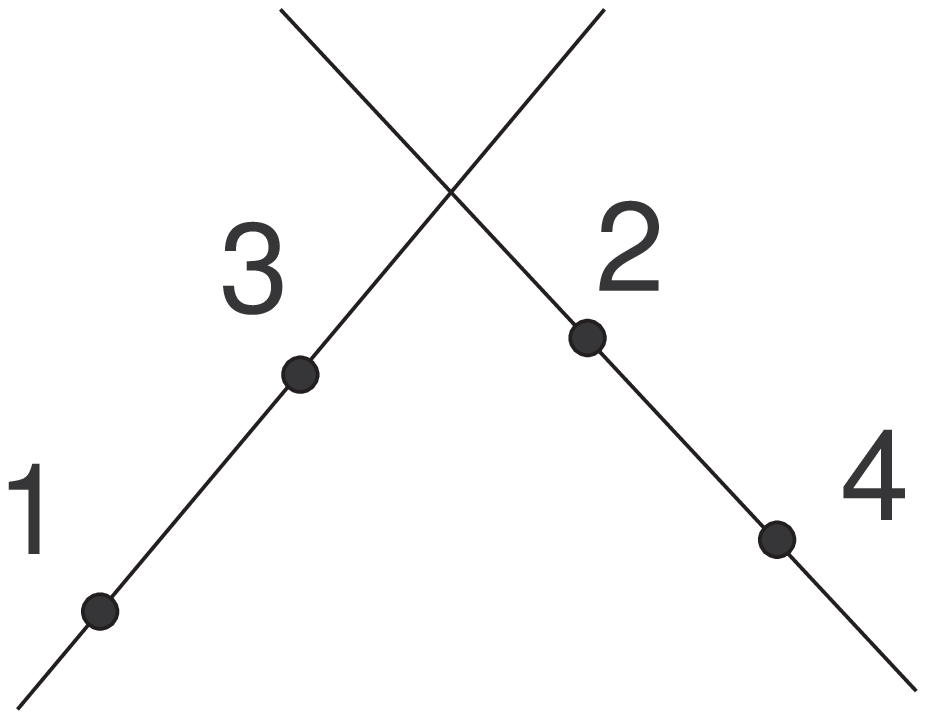}}
		& \in \overline{M}_{0,4},
\end{array}
$$
we obtain the relation:
\begin{eqnarray}
	\sum_{D_\alpha \in \varphi^{-1}(P_1)}\int_{D_\alpha} e(Ob_{\mid D_\alpha})=
	\sum_{D_\beta \in \varphi^{-1}(P_2)}\int_{D_\beta} e(Ob_{\mid D_\beta})
	\label{WDVV}
\end{eqnarray}

\noindent\textbf{Important remarks:}
\begin{enumerate}
	\item the divisors $D_\alpha$ (resp. $D_\beta$) correspond to all possible ways of distributing the remaining $n_1+n_2$ moving marks on the two twigs (Figure \ref{wdvvfig}).
	
\begin{figure}[hbt]
	\begin{center}
	\psfrag{w}{$\omega$}
\psfrag{b}{$\bar\omega$}
		\includegraphics[width=0.80\textwidth]{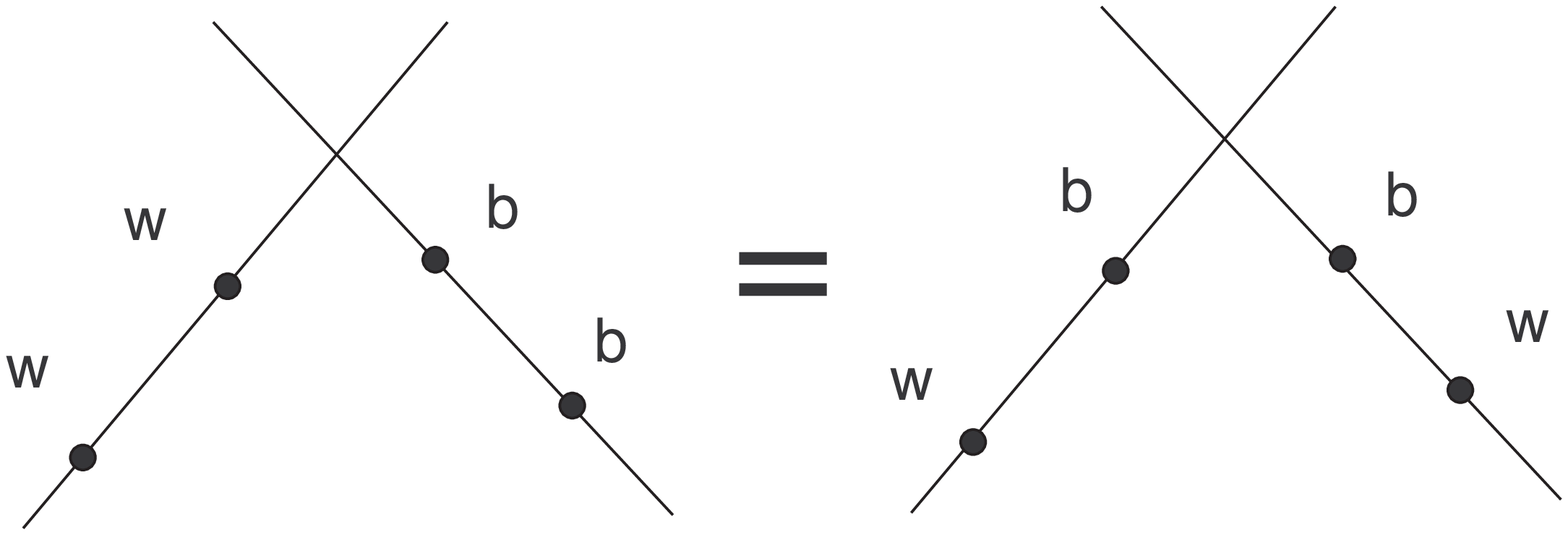}
	\end{center}
	\caption{A schematic representation of WDVV. The moving marks must be distributed on the twigs in all possible ways.}
	\label{wdvvfig}
\end{figure}
	\item Boundary divisors are (essentially)\footnote{We ignore here the gluing factor discussed in section \ref{restr}, as in this case it cancels out of the relation.} products of moduli spaces of admissible $\mathbb{Z}_3$-covers with fewer numbers of marks.
	\item The twisting at the node is determined by the monodromy condition on each twig. 
	\end{enumerate}

Next we discuss how the bundle $\mathbb{E}_\omega$ restricts to a boundary divisor $D=\mathcal{A}_1 \times\mathcal{A}_2$:
\begin{description}
	\item[Case 1:] the node is twisted. In this case
	$$
	{\mathbb{E}_\omega}_{\mid D}=\mathbb{E}^1_\omega \oplus \mathbb{E}^2_\omega,
	$$
	where we denote by $\mathbb{E}^i_\omega$ the corresponding bundle on the space $\mathcal{A}_i$.
	In this case
	$$
	e(Ob_{\mid D})=e(Ob^1)e(Ob^2).
	$$
	\item[Case 2:] the node is untwisted. Then
	$$
	{\mathbb{E}_\omega}_{\mid D}=\mathbb{E}^1_\omega \oplus \mathbb{E}^2_\omega \oplus \mathcal{O},
	$$
	where $\mathcal{O}$ is a trivial (but not equivariantly trivial!) line bundle.
	In this case
	$$
	e(Ob_{\mid D})=t_1t_2t_3\ e(Ob^1)e(Ob^2).
	$$
\end{description}

Combining these observations we obtain the following:

\begin{lemma}
WDVV gives homogeneous quadratic equations among the \GW\ invariants of $\orb$.
\end{lemma}

\noindent\textbf{Example: the smallest WDVV equation.}

Let us consider the case $n_1=n_2=0$. There are no moving points: only two divisors appear in the WDVV equation.
\begin{eqnarray}
	\langle\omega^{3}\rangle\langle\bar\omega^{3}\rangle=t_1 t_2 t_3\langle1 \omega\bar\omega\rangle^2
	\label{wdvv1}
\end{eqnarray}
Notice that this equation is consistent with our computation of three-pointed invariants in section \ref{small}. This is the only WDVV relation that features divisors with an untwisted node.
 \begin{prop}
 Let $N\geq 4$ and assume known all invariants with total number of insertions strictly less than $N$. WDVV gives a linear system of equations among all nontrivial invariants  with $N$ insertions ($N$-invariants). No invariant is directly determined by this system, but the rank of the system is one less than the number of unknowns.
 \label{wdvvrel}
 \end{prop}
 This means, once one $N$-invariant is known,  WDVV determines all other $N$-invariants inductively.
 
\begin{Prf} Consider the WDVV equation (\ref{wdvveq}) when
$
n_1+ n_2+4= N+1
$. The principal terms for this equation correspond to all moving points on the same twig. Refer to Figure \ref{wdvvfig} to analyze all possible cases:
\begin{description}
	\item [LHS, all points go left:] $\inv{n_1+3}{n_2}\langle\bar\omega^3\rangle$.
	\item [LHS, all points go right:] $\langle\omega^3\rangle\inv{n_1}{n_2+3}$.
	\item [RHS:]  principal terms have an untwisted node, and give a product of invariants with one fundamental class insertion. Since we assumed $N\geq 4$, these terms vanish. 
\end{description}
Substituting the known values for three-pointed invariants, equation (\ref{wdvveq}) reads:
\begin{eqnarray*}
	t_1t_2t_3\inv{n_1+3}{n_2}+\inv{n_1}{n_2+3}=\mbox{inductively known terms}
\end{eqnarray*}
Considering all possible values for $n_1$ and $n_2$, one obtain a matrix for the linear system:
$$
W=
\left[
\begin{array}{ccccccccc}
t_1t_2t_3 & 1 & 0 & \cdots & & & & & \\ & & & & & & & & \\
0 &t_1t_2t_3 & 1 & 0 & \cdots & & & &  \\& & & & & & & & \\
 & & & & \cdots & & & & \\& & & & & & & & \\
 & & & & \cdots& 0 & t_1t_2t_3 & 1 & 0 \\& & & & & & & & \\
  & & & & & \cdots& 0 & t_1t_2t_3 & 1  \\ 
\end{array}
\right]
$$ 
$W$ is clearly an $m-1\times m$ matrix of maximal rank, satisfying the statement of the proposition.

\end{Prf}

\section{Localization Relations}
\label{locrel}

\subsection{Maps to the trivial gerbe}
Here we obtain relations between $\mathbb{Z}_3$-Hodge integrals via localization on moduli spaces of maps of degree $1$ to a trivial gerbe (see section \ref{glocal}). These moduli spaces can be interpreted as moduli spaces of admissible covers of a parameterized $\proj$ and have been used by the second author in \cite{r:tqft}, \cite{r:dd2}, \cite{bct:gg-1}.

\subsubsection{Localization on $\tg{3k+3}{0}$}
For $k>0$, consider the auxiliary integral:
 
\begin{eqnarray*}
I_{3k+3}= \int_{\tg{3k+3}{0}}e(\E)\cup ev_1^\ast(\infty)=0
\end{eqnarray*}

$I_{3k+3}$ vanishes for dimension reasons: the degree of the integrand is $3k+2$ while the dimension of $\tg{3k+3}{0}$
is $3k+3$. We now evaluate $I_{3k+3}$ via localization and obtain relations between $\mathbb{Z}_3$-Hodge integrals. We choose to linearize the bundles according to the following table:
\begin{center}
\begin{tabular}{|l||c|c|}
\hline
weight : & over $0$ & over $\infty$ \\
\hline
\hline
${\mathcal{O}_{\proj}}$  & 0  & 0  \\
\hline${\mathcal{O}_{\proj}}$  & 0  & 0  \\
\hline
 ${\mathcal{O}_{\proj}}(-1)$ & 0 & 1 \\
\hline
\end{tabular}
\end{center}

\begin{figure}[htb]
$$
\begin{array}{clcccll}	
F_{k\times} & = & \hspace{1cm}&
	\psfrag{i}{$\infty$}
	\psfrag{0}{\vspace{.1cm}$0$}
	\psfrag{l}{\hspace{-1cm}$\adm{3k+3}{0}$}
	\psfrag{r}{\hspace{-0cm}$\omega$}
		\parbox[c]{2cm}{\includegraphics[width=.3\textwidth]{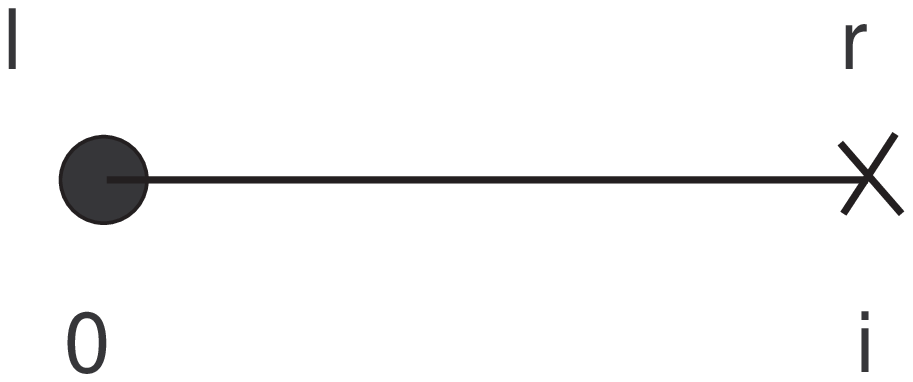}} & &\\ 
		& & & & & & \\	& & & & & & \\
F_{k_1k_2} & = & \hspace{1cm}&
	\psfrag{i}{$\infty$}
	\psfrag{0}{\vspace{.1cm}$0$}
	\psfrag{l}{\hspace{-1cm}$\adm{3k_1+1}{1}$}
	\psfrag{r}{\hspace{-.5cm}$\adm{3k_2+3}{0}$}
		\parbox[c]{2cm}{\includegraphics[width=.3\textwidth]{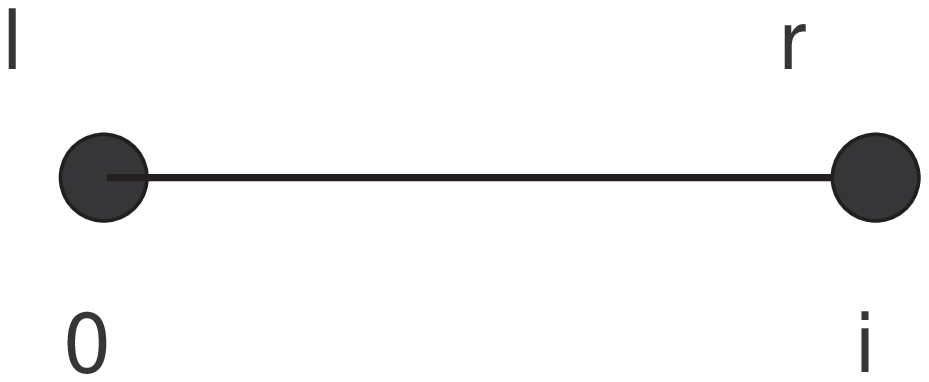}} & \hspace{2cm}&  k_1+k_2=k, & 1\leq k_1\leq k\\
		& & & & & & \\	& & & & & & \\
F_{\times k} & = & \hspace{1cm}&
	\psfrag{i}{$\infty$}
	\psfrag{0}{\vspace{.1cm}$0$}
	\psfrag{l}{\hspace{-0cm}$\omega$}
	\psfrag{r}{\hspace{-.5cm}$\adm{3k+3}{0}$}
		\parbox[c]{2cm}{\includegraphics[width=.3\textwidth]{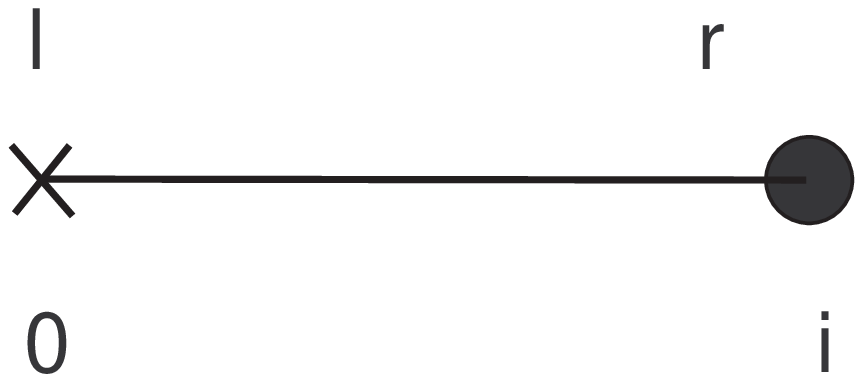}} & & \\
\end{array}
$$	
	\caption{The contributing fixed loci in the localization computation of $I_{3k+3}$.}	
	\label{loci1}
\end{figure}

This linearization and the choice of ``sending'' one point to $\infty$ force the vanishing of the contributions of many fixed loci (see \ref{restr}). The survivers are illustrated in Figure \ref{loci1} and are characterized by:
\begin{itemize}
	\item nodes over $0$ and $\infty$ are both twisted;
	\item contracted curves over $\infty$  have only $\omega$ insertions.
\end{itemize}
 
The fixed loci contributions are explicitly constructed from the data in Table \ref{ttg1}.

\begin{table}[htb]
\footnotesize
$$
\begin{array}{|c|c|c|c|c|c|}
\hline
\mbox{Locus} & \cong & \mbox{\#}  & \mbox{Edge} & \mbox{V}_0 & \mbox{V}_\infty\\
\hline
\hline
& & & & & \\
F_{k\times}: & \adm{3k+3}{0} & 1  & -\frac{2}{3}\hbar & \frac{ (-1)^k\lambda_k^2\Lambda_k(1)}{\hbar(\hbar-\psi_\omega)}& -\frac{1}{\hbar}\\
& & & & & \\
\hline
& & & & & \\
F_{k_1k_2}: & 3$\footnotesize{$\adm{3k_1+1}{1}\times\adm{3k_2+3}{0}$}$ & {3k+2}\choose{3k_1+1}  & -\frac{1}{3}\hbar & \frac{ (-1)^{k_1}\lambda_{k_1}^2\Lambda_{k_1}(1)}{\hbar(\hbar-\psi_{\bar\omega)}}& \frac{(-1)^{k_2}\lambda_{k_2}^3}{-\hbar(-\hbar-\psi_\omega)}\\
& & & & & \\
\hline
& & & & & \\
F_{\times k}: & \adm{3k+3}{0} & 3k+2  & -\frac{1}{3}\hbar & \frac{1}{\hbar}& \frac{(-1)^{k}\lambda_{k}^3}{-\hbar(-\hbar-\psi_\omega)}\\
& & & & & \\
\hline
\end{array}
$$
\normalsize
\caption{Fixed loci contributions. }
\label{ttg1}
 \end{table}

Recalling that all contributions must add to $0$, we obtain:
\begin{rel}  
\begin{eqnarray*}
		3k(\lambda_k^3)_{\adm{3k+3}{0}} &=& 2\sum_{i=1}^k(\lambda_k^2\lambda_{k-i}\psi_\omega^i)_{\adm{3k+3}{0}} -\\ \nonumber
	& & \hspace{-1.5cm}3\sum_{k_1=1}^k
	{{3k+2}\choose{3k_1+1}}\left(\sum_{j=1}^{k_1}(\lambda_{k_1}^2\lambda_{k_1-j}\psi_{\bar\omega}^{j-1})_{\adm{3k_1+1}{1}}\right)(\lambda_{k_2}^3)_{\adm{3k_2+3}{0}}
\end{eqnarray*}
\label{reltg1}
\end{rel}

\subsubsection{Localization on $\tg{3k+1}{1}$}

For $k>0$, consider:
 
\begin{eqnarray*}
I_{3k+1}= \int_{\tg{3k+1}{1}}e(\E)\cup ev_{\bar\omega}^\ast(\infty)
\end{eqnarray*}

The rank of the integrand is $3k+2$, equal to the dimension of $\tg{3k+1}{1}$. When evaluating $I_{3k+1}$, we can  linearize the three bundles arbitrarily: the result should be independent of the linearizations. For an arbitrary weight $\alpha$ we linearize the bundles according to the following table:
\begin{center}
\begin{tabular}{|l||c|c|}
\hline
weight : & over $0$ & over $\infty$ \\
\hline
\hline
${\mathcal{O}_{\proj}}$  & 0  & 0  \\
\hline${\mathcal{O}_{\proj}}$  &$\alpha$  & $\alpha$  \\
\hline
 ${\mathcal{O}_{\proj}}(-1)$ & -1 & 0 \\
\hline
\end{tabular}
\end{center}

This choice induces the vanishing of the contributions of many fixed loci. The  possibly contributing fixed loci are characterized by all nodes being  twisted (see Figure \ref{tg2fig}) and the corresponding contributions are listed in Table \ref{ttg2} .

\begin{figure}[hbt]
$$
\begin{array}{clcccll}	
F_{k\times} & = & \hspace{1cm}&
	\psfrag{i}{$\infty$}
	\psfrag{0}{\vspace{.1cm}$0$}
	\psfrag{l}{\hspace{-1cm}$\adm{3k+1}{1}$}
	\psfrag{r}{\hspace{-0cm}$\bar\omega$}
		\parbox[c]{2cm}{\includegraphics[width=.3\textwidth]{gleft.eps}} & &\\ 
		& & & & & & \\	& & & & & & \\
F_{k_1k_2} & = & \hspace{1cm}&
	\psfrag{i}{$\infty$}
	\psfrag{0}{\vspace{.1cm}$0$}
	\psfrag{l}{\hspace{-1cm}$\adm{3k_1+1}{1}$}
	\psfrag{r}{\hspace{-.5cm}$\adm{3k_2+1}{1}$}
		\parbox[c]{2cm}{\includegraphics[width=.3\textwidth]{gboth.eps}} & \hspace{2cm}&  k_1+k_2=k, & 1\leq k_1\leq k-1\\
		& & & & & & \\	& & & & & & \\
F_{\times k} & = & \hspace{1cm}&
	\psfrag{i}{$\infty$}
	\psfrag{0}{\vspace{.1cm}$0$}
	\psfrag{l}{\hspace{-0cm}$\omega$}
	\psfrag{r}{\hspace{-.5cm}$\adm{3k+1}{1}$}
		\parbox[c]{2cm}{\includegraphics[width=.3\textwidth]{gright.eps}} & & \\
		& & & & & & \\	& & & & & & \\
\tilde{F}_{k_1k_2} & = & \hspace{1cm}&
	\psfrag{i}{$\infty$}
	\psfrag{0}{\vspace{.1cm}$0$}
	\psfrag{l}{\hspace{-1cm}$\adm{3k_1+3}{0}$}
	\psfrag{r}{\hspace{-.5cm}$\adm{3k_2-1}{2}$}
		\parbox[c]{2cm}{\includegraphics[width=.3\textwidth]{gboth.eps}} & \hspace{2cm}&  k_1+k_2=k, & 0\leq k_1\leq k-1\\
		& & & & & & \\	& & & & & & \\
\end{array}
$$	
	\caption{The contributing fixed loci in the localization computation of $I_{3k+1}$.}	
	\label{tg2fig}
\end{figure}	
\begin{table}[!h]
\footnotesize
$$
\begin{array}{|c|c|c|c|c|c|}
\hline
\mbox{Locus} & \cong & \mbox{\#}  & \mbox{Edge} & \mbox{V}_0 & \mbox{V}_\infty\\
\hline
\hline
& & & & & \\
F_{k\times}: & \adm{3k+1}{1} & 1  & -\frac{1}{3}\hbar & \frac{ (-1)^k\lambda_k\Lambda_k(-\alpha)\Lambda_k(1)}{\hbar(\hbar-\psi_{\bar\omega)}}& -\frac{1}{\hbar}\\
& & & & & \\
\hline
& & & & & \\
F_{k_1k_2}: & 3$\footnotesize{$\adm{3k_1+1}{1}\times\adm{3k_2+1}{1}$}$ & {3k+1}\choose{3k_1+1}  & -\frac{1}{3}\hbar & \frac{ (-1)^{k_1}\lambda_{k_1}\Lambda_{k_1}(-\alpha)\Lambda_{k_1}(1)}{\hbar(\hbar-\psi_{\bar\omega)}}& \frac{(-1)^{k_2}\lambda_{k_2}^2\Lambda_{k_2}(-\alpha)}{-\hbar(-\hbar-\psi_\omega)}\\
& & & & & \\
\hline
& & & & & \\
F_{\times k}: & \adm{3k+1}{1} & 3k+1  & -\frac{1}{3}\hbar & \frac{1}{\hbar}& \frac{(-1)^{k}\lambda_{k}^2\Lambda_{k_1}(-\alpha)}{-\hbar(-\hbar-\psi_\omega)}\\
& & & & & \\
\hline
& & & & & \\
\tilde{F}_{k_1k_2}: & 3$\footnotesize{$\adm{3k_1+3}{0}\times\adm{3k_2-1}{2}$}$ & {3k+1}\choose{3k_1+2}  & -\frac{2}{3}\hbar & \frac{(-1)^{k_1}\lambda_{k_1}\Lambda_{k_1}(-\alpha)\Lambda_{k_1}(1)}{\hbar(\hbar-\psi_\omega)}& \frac{(-1)^{k_2}\lambda_{k_2}^2\Lambda_{k_2}(-\alpha)}{-\hbar(-\hbar-\psi_{\bar\omega})}\\
& & & & & \\
\hline
\end{array}
$$
\normalsize
\caption{Fixed loci contributions. }
\label{ttg2}
 \end{table}

For $\alpha=0$, $F_{k\times}$ is the only contributing fixed locus and
\begin{eqnarray}
	I_{3k+1}=\frac{(-1)^k}{3}\sum_{i=1}^k (\lambda_k^2\lambda_{k-i}\psi_{\bar\omega}^{i-1})_{\adm{3k+1}{1}}
\label{alpha0}
\end{eqnarray}

Subtracting (\ref{alpha0}) from the evaluation of  $I_{3k+1}$ for a general value of $\alpha$ we obtain a polynomial in $\alpha$. All of its coefficients must vanish thus giving relations among $\mathbb{Z}_3$-Hodge integrals. We focus on the vanishing of the linear coefficient:  

\begin{rel}  

\begin{eqnarray*}
3k(\lambda_k^2\lambda_{k-1})_{\adm{3k+1}{1}}-(\lambda_k\lambda_{k-1}^2\psi_{\bar\omega})_{\adm{3k+1}{1}}=	 \sum_{i=2}^k\hi{k}{k-1}{k-i}{\bar\omega}{i-2}{3k+1}{1}-
	\nonumber\\	3\sum_{k_1=1}^{k-1}{{3k+1}\choose{3k_1+1}}\left(\sum_{i=1}^{k_1}\hii{{k_1}}{{k_1}-i}{\bar\omega}{i-1}{3{k_1}+1}{1}\right)\hiinop{{k_2}}{{k_2}-1}{3{k_2}+1}{1}
\end{eqnarray*}
\label{reltg2}
\end{rel}
\textbf{Remarks:}
\begin{enumerate}
	\item We have chosen to isolate two terms that will play the role of principal parts in our inductive procedure for computing all invariants of $\orb$.
	\item Notice that by choosing to look only at the linear part in $\alpha$ we have no contributions from the $\tilde{F}_{k_1k_2}$ loci. 
\end{enumerate}
	
\subsection{Relations from maps to $G_1$}

Here we apply localization to moduli spaces of maps to the first $\mathbb{Z}_3$-gerbe over $\proj$. Even though there is no ``non-stacky'' interpretation for the general map in these moduli spaces, the fixed loci are (essentially) products of spaces admissible covers. Thus we extract relations between $\mathbb{Z}_3$-Hodge integrals.  

\subsubsection{Localization on $\fg{3k+1}{0}$}
For $k>0$, consider the auxiliary integral:
 
\begin{eqnarray*}
J_{3k+1}= \int_{\fg{3k+1}{0}}e\left(\EE\right) = 0
\end{eqnarray*}

$J_{3k+1}$ vanishes for dimension reasons: the rank of the integrand is $3k$, the dimension of $\fg{3k+1}{0}$
is $3k+1$.  We linearize the bundles:
\begin{center}
\begin{tabular}{|l||c|c|}
\hline
weight : & over $0$ & over $\infty$ \\
\hline
\hline
${\mathcal{O}_{G_1}(-2/3)}$  & 0  & 2/3  \\
\hline${\mathcal{O}_{G_1}(-2/3)}$  & -2/3  & 0  \\
\hline
 ${\mathcal{O}_{G_1}}(-2/3)$ & 0 & 2/3\\
\hline
\end{tabular}
\end{center}

\begin{figure}[hbt]
$$
\begin{array}{clcccll}	
F_{k\circ} & = & \hspace{1cm}&
	\psfrag{i}{$\infty$}
	\psfrag{0}{\vspace{.1cm}$0$}
	\psfrag{l}{\hspace{-1cm}$\adm{3k+1}{1}$}
	\psfrag{r}{\hspace{-0cm}$1$}
		\parbox[c]{2cm}{\includegraphics[width=.3\textwidth]{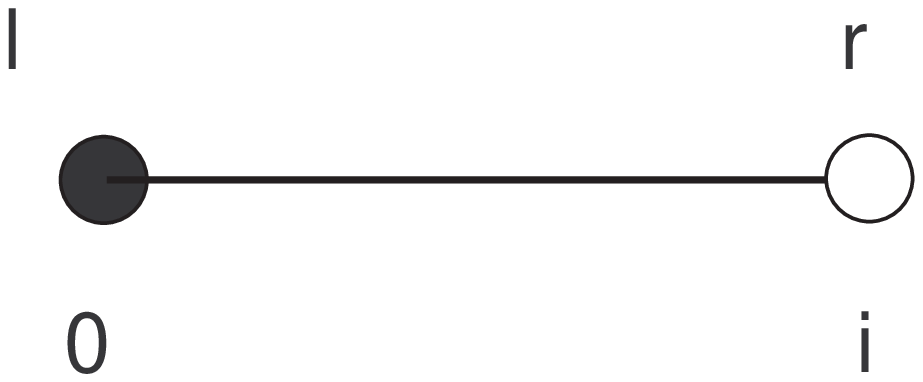}} & &\\ & & & & & & \\
		F_{\circ k} & = & \hspace{1cm}&
	\psfrag{i}{$\infty$}
	\psfrag{0}{\vspace{.1cm}$0$}
	\psfrag{l}{\hspace{-0cm}$1$}
	\psfrag{r}{\hspace{-.5cm}$\adm{3k+1}{1}$}
		\parbox[c]{2cm}{\includegraphics[width=.3\textwidth]{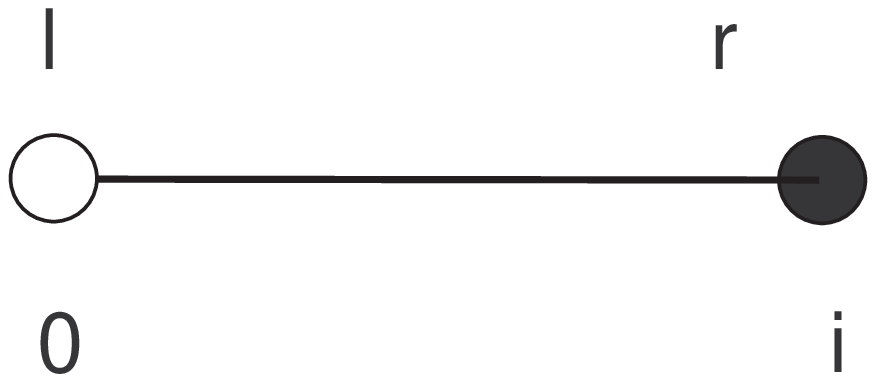}} & &\\  & & & & & & \\
F_{k_1k_2} & = & \hspace{1cm}&
	\psfrag{i}{$\infty$}
	\psfrag{0}{\vspace{.1cm}$0$}
	\psfrag{l}{\hspace{-1cm}$\adm{3k_1+3}{0}$}
	\psfrag{r}{\hspace{-.5cm}$\adm{3k_2+3}{0}$}
		\parbox[c]{2cm}{\includegraphics[width=.3\textwidth]{gboth.eps}} & \hspace{2cm}&  k_1+k_2=k-1, & 0\leq k_1\leq k-1\\
\end{array}
$$	
	\caption{The contributing fixed loci in the localization computation of $J_{3k+1}$.}	
	\label{fg1fig}
\end{figure}

Figure \ref{fg1fig} and Table \ref{fgt1} illustrate the nonvanishing fixed loci and their contributions.
 
\begin{table}[hbt]
\footnotesize
$$
\begin{array}{|c|c|c|c|c|c|}
\hline
\mbox{Locus} & \cong & \mbox{\#}  & \mbox{Edge} & \mbox{V}_0 & \mbox{V}_\infty\\
\hline
\hline
& & & & & \\
F_{k\circ}: & \adm{3k+1}{1} & 1  & 1 & \frac{ (-1)^k\lambda_k^2\Lambda_k(2/3)}{\hbar(\hbar-\psi_{\bar\omega})}& 1\\
& & & & & \\
\hline
& & & & & \\
F_{\circ k}: & \adm{3k+1}{0} & 1  & 1 & 1 & \frac{(-1)^{k}\lambda_{k}\Lambda_{k}(-2/3)^2}{-\hbar(-\hbar-\psi_{\bar\omega})}\\
& & & & & \\
\hline
& & & & & \\
F_{k_1k_2}: & 3$\footnotesize{$\adm{3k_1+3}{0}\times\adm{3k_2+3}{0}$}$ & {3k+1}\choose{3k_1+2}  & -\frac{1}{27}\hbar^3 & \frac{ (-1)^{k_1}\lambda_{k_1}^2\Lambda_{k_1}(2/3)}{\hbar(\hbar-\psi_\omega)}& \frac{(-1)^{k_2}\lambda_{k_2}\Lambda_{k_2}(-2/3)^2}{-\hbar(-\hbar-\psi_\omega)}\\
& & & & & \\
\hline
\end{array}
$$
\normalsize
\caption{Fixed loci contributions.}
\label{fgt1}
\end{table}

Recalling that all contributions must add to $0$, we obtain:
\begin{rel}  

\begin{eqnarray*}
		6(\lambda_k^2\lambda_{k-1})_{\adm{3k+1}{1}}+4(\lambda_k\lambda_{k-1}^2)_{\adm{3k+1}{1}} &=&   9\left[\sum_{i=2}^k\left(\frac{2}{3}\right)^i(\lambda_k^2\lambda_{k-i}\psi^{i-1}_{\bar\omega})_{\adm{3k+1}{1}}- \right.\nonumber\\
& &		
\hspace{-3cm}\left. 		\sum_\star\left(\frac{2}{3}\right)^{i+j}\hi{k}{k-i}{k-j}{\bar{\omega}}{i+j-1}{3k+1}{1}
		\right]+   \nonumber\\
	& &\hspace{-3cm}  \sum_{k_1=0}^{k-1}
{{3k+1}\choose{3k_1+2}}\left(\sum_{i=0}^{k_1}\left(\frac{2}{3}\right)^i(\lambda_{k_1}^2\lambda_{k_1-i}\psi_{\omega}^{j})_{\adm{3k_1+3}{0}}\right)\cdot \nonumber \\
 & & \hspace{-3cm}\cdot
\left(\sum_{i,j=0}^{k_1}\left(\frac{2}{3}\right)^{i+j}\hi{k_2}{k_2-i}{k_2-j}{\omega}{i+j}{3k_2+3}{0}\right)
\end{eqnarray*}
\label{relfg1}
\end{rel}	

The symbol $\star$ stands for: $0\leq i,j\leq k$, $i+j\geq2$, $(i,j)\not=(1,1)$.
 
\subsubsection{Localization on $\fg{3k-1}{1}$}

For $k>0$, consider:
 
\begin{eqnarray*}
J_{3k-1}= \int_{\fg{3k-1}{1}}e\left(\EE\right)
\end{eqnarray*}

The rank of the integrand is $3k$, equal to the dimension of $\fg{3k-1}{1}$. The evaluation of the integral is independent of the linearization of the bundles.
For an arbitrary weight $\alpha$ we choose:
\begin{center}
\begin{tabular}{|l||c|c|}
\hline
weight : & over $0$ & over $\infty$ \\
\hline
\hline
${\mathcal{O}_{G_1}}(-2/3)$  & 0  & 2/3  \\
\hline${\mathcal{O}_{G_1}}(-2/3)$  &$-2/3$  & $0$  \\
\hline
 ${\mathcal{O}_{G_1}}(-2/3)$ & $\alpha$ & $\alpha+2/3$ \\
\hline
\end{tabular}
\end{center}

This choice induces the vanishing of the contributions of many fixed loci. The  possibly contributing fixed loci are those with no untwisted nodes (see Figure \ref{fg2fig}, Table \ref{fgt2}).
 
\begin{figure}[hbt]
$$
\begin{array}{clcccll}	
F_{k\circ} & = & \hspace{1cm}&
	\psfrag{i}{$\infty$}
	\psfrag{0}{\vspace{.1cm}$0$}
	\psfrag{l}{\hspace{-1cm}$\adm{3k-1}{2}$}
	\psfrag{r}{\hspace{-0cm}$1$}
		\parbox[c]{2cm}{\includegraphics[width=.3\textwidth]{gleftc.eps}} & &\\ 
				& & & & & & \\	& & & & & & \\
F_{\circ k} & = & \hspace{1cm}&
	\psfrag{i}{$\infty$}
	\psfrag{0}{\vspace{.1cm}$0$}
	\psfrag{l}{\hspace{-0cm}$1$}
	\psfrag{r}{\hspace{-.5cm}$\adm{3k-1}{2}$}
		\parbox[c]{2cm}{\includegraphics[width=.3\textwidth]{grightc.eps}} & & \\
		& & & & & & \\	& & & & & & \\
F_{k_1k_2} & = & \hspace{1cm}&
	\psfrag{i}{$\infty$}
	\psfrag{0}{\vspace{.1cm}$0$}
	\psfrag{l}{\hspace{-1cm}$\adm{3k_1+1}{1}$}
	\psfrag{r}{\hspace{-.5cm}$\adm{3k_2+3}{0}$}
		\parbox[c]{2cm}{\includegraphics[width=.3\textwidth]{gboth.eps}} & \hspace{2cm}&  k_1+k_2=k-1, & 0 \leq k_1\leq k-1\\
		& & & & & & \\	& & & & & & \\
\tilde{F}_{k_1k_2} & = & \hspace{1cm}&
	\psfrag{i}{$\infty$}
	\psfrag{0}{\vspace{.1cm}$0$}
	\psfrag{l}{\hspace{-1cm}$\adm{3k_1+3}{0}$}
	\psfrag{r}{\hspace{-.5cm}$\adm{3k_2+1}{1}$}
		\parbox[c]{2cm}{\includegraphics[width=.3\textwidth]{gboth.eps}} & \hspace{2cm}&  k_1+k_2=k-1, & 0\leq k_1\leq k-1\\
		& & & & & & \\	& & & & & & \\
\end{array}
$$	
	\caption{The contributing fixed loci in the localization computation of $J_{3k-1}$.}	
	\label{fg2fig}
\end{figure}	
	
\begin{table}[hbt]
\footnotesize
$$
\hspace{-3cm}\begin{array}{|c|c|c|c|c|c|}
\hline
\mbox{Locus} & \cong & \mbox{\#}  & \mbox{Edge} & \mbox{V}_0 & \mbox{V}_\infty\\
\hline
\hline
& & & & & \\
F_{k\circ}: & \adm{3k-1}{2} & 1  & 1 & \frac{ (-1)^k\lambda_k\Lambda_k(-\alpha)\Lambda_k(2/3)}{\hbar(\hbar-\psi_{\bar\omega)}}& 1\\
& & & & & \\
\hline
& & & & & \\
F_{\circ k}: & \adm{3k-1}{2} & 1  & 1 & 1& \frac{(-1)^{k}\lambda_{k}\Lambda_{k}(-\alpha-2/3)\Lambda_{k}(-2/3)}{-\hbar(-\hbar-\psi_{\bar\omega)}}\\
& & & & & \\
\hline
& & & & & \\
F_{k_1k_2}: & 3$\footnotesize{$\adm{3k_1+1}{1}\times\adm{3k_2+3}{0}$}$ & {3k-1}\choose{3k_1}  &  -\frac{1}{9}\left(\alpha+\frac{1}{3}\right)\hbar^3  & \frac{ (-1)^{k_1}\lambda_{k_1}\Lambda_{k_1}(-\alpha)\Lambda_{k_1}(2/3)}{\hbar(\hbar-\psi_\omega)}& \frac{(-1)^{k_2}\lambda_{k_2}^2\Lambda_{k_2}(-\alpha-2/3)\Lambda_{k_2}(-2/3)}{-\hbar(-\hbar-\psi_\omega)}\\
& & & & & \\
\hline
& & & & & \\
\tilde{F}_{k_1k_2}: & 3$\footnotesize{$\adm{3k_1+3}{0}\times\adm{3k_2+1}{1}$}$ & {3k-1}\choose{3k_1+2}  & -\frac{1}{9}\left(\alpha+\frac{1}{3}\right)\hbar^3 & \frac{ (-1)^{k_1}\lambda_{k_1}\Lambda_{k_1}(-\alpha)\Lambda_{k_1}(2/3)}{\hbar(\hbar-\psi_\omega)}& \frac{(-1)^{k_2}\lambda_{k_2}^2\Lambda_{k_2}(-\alpha-2/3)\Lambda_{k_2}(-2/3)}{-\hbar(-\hbar-\psi_\omega)}\\
\hline
\end{array}
$$
\normalsize
\caption{Fixed loci contributions.}
\label{fgt2}
\end{table}
\textbf{Note:} the contribution of the degenerate locus $F_{0,k}$ (resp. $\tilde{F}_{k-1,0}$) can be read from  Table \ref{fgt2} by adopting the convention that the contribution of $V_0$ (resp. $V_\infty$) is defined to be $\frac{1}{3}$ (resp. $-\frac{1}{3}$).

Recalling that $J_{3k-1}(\alpha) - J_{3k-1}(0)=0$ as a polynomial in $\alpha$, we obtain a relation from the vanishing of the first degree coefficient.
\begin{rel}
\begin{eqnarray*}
	\frac{4}{3}(\lambda_{k}^2\lambda_{k-2})_{\adm{3k-1}{2}} &= &  \sum_{i=2}^k \hi{k}{k-1}{k-i}{{\bar\omega}}{i-1}{3k-1}{2}- \nonumber\\
	& &  \hspace{-1cm}\sum_{\star}i \left(\frac{2}{3}\right)^{i+j-1}\hi{k}{k-i}{k-j}{{\bar\omega}}{i+j-1}{3k-1}{2}  + \nonumber \\
	& &  \hspace{-1cm} \ \mbox{terms on strictly ``smaller'' moduli spaces}
\end{eqnarray*}
\label{relfg2}
\end{rel}
\textbf{Remarks:}
\begin{enumerate}
	\item here $\star$ means $0\leq i,j \leq k, i+j>2$.
	\item we choose not to record the full relation here simply because it is longer than the previous ones. It is not however more (computationally) complex. 
\end{enumerate}
\subsection{Removing Descendant Insertions}
\label{rd}
In this section we give a series of recursions that express any $\mathbb{Z}_3$-Hodge integral of the form
\begin{eqnarray}
	\int_{\adm{n_1}{n_2\leq 2}}\lambda_{k}\lambda_{k-i}\lambda_{k-j}\psi^{l}
	\label{psi}
\end{eqnarray}
(where $k$ is the rank of the $\mathbb{E}_\omega$ in question and $l$ is strictly positive) in terms of integrals on strictly smaller moduli spaces. The strategy is the same as in section \ref{locrel}. Since we feel we have already provided a sufficient amount of detailed localization computations, here we only state the  vanishing auxiliary integrals, and we expand only one example that we specifically need in the proof of our main result. 
\begin{notation}
In the following paragraph we adopt the notation:
\begin{itemize}
	\item  $\lambda_i$ to mean $c_i((R^1\pi_\ast f^\ast(\mathcal{O}_{\proj}))^\vee)$. We also assume $\mathcal{O}_{\proj}$ linearized with weights $(0,0)$.
	\item $ev_\omega$ to indicate an evaluation map corresponding to an $\omega$ point (likewise for $\bar\omega$).	
\end{itemize}
\end{notation}
\begin{description}
\item[(a) Removing $\psi_\omega$'s from $\adm{3k+3}{0}$.] To compute $\hi{k}{k-i}{k-j}{\omega}{i+j}{3k+3}{0}$ we use the auxiliary integral:
\begin{eqnarray*}
\int_{\tg{3k+3}{0}}\lambda_{k}\lambda_{k-i}\lambda_{k-j} \cup ev_\omega(0)\cup ev_\omega(0)\cup ev_\omega(\infty) =0.
\end{eqnarray*}
\item[(b) Removing $\psi_{\bar\omega}$'s from $\adm{3k+1}{1}$.] To compute $\hi{k}{k-i}{k-j}{{\bar\omega}}{i+j-1}{3k+1}{1}$:
\begin{eqnarray*}
\int_{\tg{3k+1}{1}}\lambda_{k}\lambda_{k-i}\lambda_{k-j} \cup ev_\omega(0)\cup ev_\omega(0)\cup ev_{\bar\omega}(\infty) =0.
\end{eqnarray*}
\item[(c) Removing $\psi_{\omega}$'s from $\adm{3k+1}{1}$.] To compute $\hi{k}{k-i}{k-j}{\omega}{i+j-1}{3k+1}{1}$:
\begin{eqnarray*}
\int_{\tg{3k+1}{1}}\lambda_{k}\lambda_{k-i}\lambda_{k-j} \cup ev_\omega(0)\cup ev_{\bar\omega}(0)\cup ev_\omega(\infty) =0.
\end{eqnarray*}
\item[(d) Removing $\psi_{\bar\omega}$'s from $\adm{3k-1}{2}$.] To compute $\hi{k}{k-i}{k-j}{\omega}{i+j-2}{3k-1}{2}$:
\begin{eqnarray*}
\int_{\tg{3k-1}{2}}\lambda_{k}\lambda_{k-i}\lambda_{k-j} \cup ev_\omega(0)\cup ev_{\bar\omega}(0)\cup ev_{\bar\omega}(\infty) =0.
\end{eqnarray*}
\end{description}
\textbf{Example:} to illustrate how these recursions work  we analyze the case:
\begin{eqnarray*}
\int_{\tg{3k+1}{1}}\lambda_{k}\lambda_{k-1}\lambda_{k-1} \cup ev_\omega(0)\cup ev_{\bar\omega}(0)\cup ev_{\bar\omega}(\infty) =0	
\end{eqnarray*}
Since we require two twisted points to ``go to $0$'', only components that have a node over $0$ can contribute. Further, we must have either a node or a twisted point at $\infty$. The nontrivial contributions are illustrated in the following table:
\footnotesize
$$
\begin{array}{|c|c|c|c|c|}
\hline
\mbox{Locus} & \cong & \mbox{\#}  &  \mbox{V}_0 & \mbox{V}_\infty\\
\hline
\hline
& & & &  \\
F_{k\times} & \adm{3k+1}{1} & 1  & \frac{\lambda_k\lambda^2_{k-1}}{\hbar(\hbar-\psi_{\bar\omega)}}& -\frac{1}{\hbar}\\
& & & &  \\
\hline
& & & &  \\
\begin{array}{c}
F_{k_1k_2}\\
(1 \leq k_1\leq k-1)
\end{array} & 3$\footnotesize{$\adm{3k_1+1}{1}\times\adm{3k_2+1}{1}$}$ & {3k-1}\choose{3k_1-1}   & \frac{\lambda_{k_1}^2\lambda_{k_1-1} }{\hbar(\hbar-\psi_{\bar\omega)}}& \frac{\lambda_{k_2}^2\lambda_{{k_2}-1}}{-\hbar(-\hbar-\psi_\omega)}\\
& & & &  \\
\hline
& & & &  \\
\begin{array}{c}
\widetilde{F}_{k_1k_2}\\
(0 \leq k_1\leq k-1)
\end{array} & 3$\footnotesize{$\adm{3k_1+3}{0}\times\adm{3k_2-1}{2}$}$ & {3k-1}\choose{3k_1}& \frac{\lambda_{k_1}^3}{\hbar(\hbar-\psi_\omega)}& \frac{\lambda_{k_1}\lambda_{k_1-1}^2}{-\hbar(-\hbar-\psi_{\bar\omega})}\\
& & & &  \\
\hline
\end{array}
$$
\normalsize

\begin{rel}
\begin{eqnarray*}
	(\lambda_{k}\lambda_{k-1}^2\psi_{\bar\omega})_{\adm{3k+1}{1}} 
	-(\lambda_{k}^3)_{\adm{3k+3}{0}} (\lambda_{k}\lambda^2_{{k}-1})_{\adm{3k-1}{2}}& = & \nonumber \\ 
	 & & \hspace{-9.5cm} \sum_{k_1=1}^{k-1}(\lambda_{k_1}^2\lambda_{k_1-1})_{\adm{3k_1+1}{1}} (\lambda_{k_2}^2\lambda_{{k_2}-1})_{\adm{3k_2+1}{1}}+
	\sum_{k_1=1}^{k-1}(\lambda_{k_1}^3)_{\adm{3k_1+3}{0}} (\lambda_{k_2}\lambda^2_{{k_2}-1})_{\adm{3k_2-1}{2}} \nonumber
\end{eqnarray*}
\label{relrempsi}
\end{rel}

Iterated use of these relations yield the following simple reconstruction result.

\begin{prop}
Descendant $\mathbb{Z}_3$-Hodge integrals of the form (\ref{psi}) can be reconstructed from non-descendant $\mathbb{Z}_3$-Hodge integrals on strictly smaller moduli spaces.
\end{prop}
\section{Proof of Main Result}
\label{proof}
In these section we prove that the localization computations of section \ref{locrel} together with WDVV provide inductive recursions that allow one to (effectively) compute any equivariant \GW\ of $\orb$. The initial data required are the three pointed invariants computed in section \ref{small}.

\noindent\textbf{Reduction: minimizing the number of $\bar\omega$ insertions to consider.} By Proposition \ref{wdvvrel}, if invariants with  at most two $\bar\omega$ insertions are known, WDVV suffices to determine all other invariants. We therefore restrict our attention to such invariants. We have four classes of $\mathbb{Z}_3$-Hodge integrals to determine:
\begin{description}
	\item[on $\adm{3k+3}{0}:$] $$ \lambda_k^3;$$
	\item[on $\adm{3k+1}{1}:$] $$ \lambda_k^2\lambda_{k-1};$$
	\item[on $\adm{3k-1}{2}:$] $$ \lambda_k^2\lambda_{k-2} \mbox{and}\ \lambda_k\lambda_{k-1}^2 .$$
\end{description}

\noindent\textbf{Induction:}

\textbf{Step 1.} Assume known all invariants with strictly less than $3k+3$ insertions. Then relation \ref{reltg1} expresses $(\lambda_k^3)_{\adm{3k+3}{0}}$ in terms of (products of) strictly smaller invariants and descendant invariants on $\adm{3k+3}{0}$. These are reduced to smaller invariants by applying the recursions in section \ref{rd} We therefore know $(\lambda_k^3)_{\adm{3k+3}{0}}$.

\textbf{Step 2.} Relation \ref{relfg2} computes $(\lambda_k^2\lambda_{k-2})_{\adm{3k-1}{2}}$ in terms of known quantities.

\textbf{Step 3.} We now observe relations \ref{reltg2}, \ref{relfg1} and \ref{relrempsi}. The Hodge integrals that are not
already known after the first two steps in our induction are $(\lambda_k^2\lambda_{k-1})_{\adm{3k+1}{1}}$,
$(\lambda_k\lambda^2_{k-1}\psi_{\bar\omega})_{\adm{3k+1}{1}}$ and $(\lambda_k \lambda^2_{k-1})_{\adm{3k-1}{2}}$.
We have a linear system of three equation in three unknowns. It is immediate from our presentation of the recursions to see that it is an invertible system. We therefore know all invariants with strictly less than $3(k+1) +3$ insertions and can start over from step 1 again.

\qed

\section{PDE's controlling the \GW\ theory of $\orb$.}
\label{gf}

Generating functions are a very efficient method to package information about systems of numbers with a rich combinatorial structure. The idea is simple: the numbers are organized to be the coefficients of some formal power series, and the relations among the numbers described in terms of differential equations among these power series. In this section we present our recursions in compact generating function form.

\subsection{WDVV}

All WDVV relations of section \ref{WDVVrel} are contained in a unique homogeneous quadratic PDE on the \GW\ potential of $\orb$. Define:
$$
\mathcal{F}(x_0,x_1,x_2):= \sum_{n_0,n_1,n_2} \langle 1^{n_0} \omega^{n_1}\bar\omega^{n_2}\rangle \frac{x_0^{n_0}}{n_0!}\frac{x_1^{n_1}}{n_1!}\frac{x_2^{n_2}}{n_2!}
$$

Then, WDVV becomes:

\begin{eqnarray}
	\mathcal{F}_{x_1x_1x_1}\mathcal{F}_{x_2x_2x_2}-\mathcal{F}_{x_1x_1x_2}\mathcal{F}_{x_1x_2x_2}=t_1t_2t_3\mathcal{F}_{x_0x_1x_2}^2\left(=\frac{1}{9}t_1t_2t_3\right).
\end{eqnarray}

\subsection{Localization Relations}

Our localization relations are best expressed in terms of  appropriate generating functions for \zhh:

\begin{eqnarray}\label{gfhi}
	\mathcal{L}^\omega(x,y;u,v):= \sum_{m,n,i,j}\hi{top}{top-i}{top-j}{\omega}{i+j-n}{m}{n} \frac{x^m}{m!}\frac{y^n}{n!} u^i v^j
\end{eqnarray}
Similarly, define $\mathcal{L}^{\bar\omega}$ by replacing $\psi_\omega$ with $\psi_{\bar\omega}$ in (\ref{gfhi}). The localization relations translate to the following PDE's. 

\begin{description}
\item[Relation 1:]
$$
2\mathcal{L}^\omega_x (-x,0;1,0)= 3\mathcal{L}_y^{\bar\omega}(-x,0;1,0) \mathcal{L}^\omega_{xx}(-x,0;0,0)
$$
\item[Relation 2:]
$$\hspace{-1.5cm}
2\mathcal{L}^\omega_x (-x,0;-u,1)\mathcal{L}^{\bar\omega}_{yy}(-x,0;u,0)-\mathcal{L}^{\bar\omega}_y(-x,0;-u,1) \mathcal{L}^{\omega}_{xy}(-x,0;u,0)+\frac{1}{3}\mathcal{L}^{\bar\omega}_y(-x,0;0,1)=0
$$
\item[Relation 3:]
$$
\mathcal{L}^{\bar\omega}_y\left(-x,0;\frac{2}{3},\frac{2}{3}\right)-\mathcal{L}^{\bar\omega}_y\left(-x,0;\frac{2}{3},0\right)=\frac{1}{9}\mathcal{L}^\omega_x\left(-x,0;\frac{2}{3},\frac{2}{3}\right)\mathcal{L}^\omega_x\left(-x,0;\frac{2}{3},0\right)
$$
\item[Relation 4:]
$$\hspace{-1.5cm}
\mathcal{L}^{\bar\omega}_{yy}\left(-x,0;\frac{2}{3},v+\frac{2}{3}\right)+\mathcal{L}^{\bar\omega}_{yy}\left(-x,0;\frac{2}{3},-v\right)-\frac{1}{3}\left(v+\frac{1}{3}\right)\left[\mathcal{L}^\omega_x\left(-x,0;\frac{2}{3},-v\right) \mathcal{L}^{\omega}_{xy}\left(-x,0;\frac{2}{3},v+\frac{2}{3}\right)-   \right. $$

$$\left.-\mathcal{L}^{\omega}_{xy}\left(-x,0;\frac{2}{3},-v\right)\mathcal{L}^\omega_x\left(-x,0;\frac{2}{3},v+\frac{2}{3}\right)\right]=2\mathcal{L}^{\bar\omega}_{yy}\left(-x,0;\frac{2}{3},\frac{1}{3}\right)
$$

\item[$\psi$ Removal - (a):]
$$ \mathcal{L}^\omega_{xxx}(x,0;u,v)\mathcal{L}^{\bar\omega}_{xy} (x,0;-u,-v)= \frac{1}{9}+\mathcal{L}^{\bar\omega}_{xxy} (x,0;u,v) \mathcal{L}^\omega_{xx}(x,0;-u,-v)$$
\item[$\psi$ Removal - (b):]
$$\mathcal{L}^{\bar\omega}_{xxy} (x,0;u,v)\mathcal{L}^\omega_{xy} (x,0;-u,-v)=\mathcal{L}^\omega_{xxx}(x,0;u,v)\mathcal{L}^{\bar\omega}_{yy}(x,0;-u,-v)$$
\item[$\psi$ Removal - (c):]
$$\mathcal{L}^{\omega}_{xxy} (x,0;u,v)\mathcal{L}^{\bar\omega}_{xy} (x,0;-u,-v)=\mathcal{L}^{\bar\omega}_{xyy}(x,0;u,v)\mathcal{L}^{\omega}_{xx}(x,0;-u,-v)$$
\item[$\psi$ Removal - (d):]
$$\mathcal{L}^{\bar\omega}_{xyy}(x,0;u,v)\mathcal{L}^\omega_{xy} (x,0;-u,-v)=-\frac{1}{9}xuv +\mathcal{L}^\omega_{xxy} (x,0;u,v)\mathcal{L}^{\bar\omega}_{yy}(x,0;-u,-v) $$
\end{description}
\subsection{$\mathcal{L}$ and \GW\ invariants}

In section \ref{oizhh} we expressed \GW\ invariants of $\orb$ in terms of \zhh. Our localization relations give a system of recursions between invariants with at most $2$ $\bar\omega$ points. In terms of Hodge integrals, this yields the significant simplification that all such invariants contain at least one $\lambda_{top}$ class. These invariants can therefore be easily related to our $\mathcal{L}$ functions.\footnote{The superscript here is unnecessary since we are dealing only with primary invariants.}  Precisely, we have:
\begin{description}
	\item[Invariants with $0$ $\bar\omega$ points:]
	$$
	\mathcal{F}(0,x_1,0) =\mathcal{L}(x_1,0;0,0)
	$$
		\item[Invariants with $1$ $\bar\omega$ point:]
	$$
	\mathcal{F}_{x_2}(0,x_1,0) =(t_1+t_2+t_3)\mathcal{L}_{y u}(x_1,0;0,0)
	$$
		\item[Invariants with $2$ $\bar\omega$ points:]
	$$\hspace{-1cm}
	\mathcal{F}_{x_2x_2}(0,x_1,0) =\frac{(t_1^2+t_2^2+t_3^2)}{2}\mathcal{L}_{yy uu}(x_1,0;0,0)+ (t_1t_2+t_1t_3+t_2t_3)\mathcal{L}_{yy uv}(x_1,0;0,0)
	$$
\end{description}
\begin{Rem}
The expert eye will notice that the generating functions $\mathcal{L}^\omega,\mathcal{L}^{\bar\omega}$ are very closely related to Givental's (equivariant) $J$ function:
$$
J(x,y;t_1,t_2,t_3):= \sum_{n} \frac{1}{n!}\langle x\omega+y\bar\omega,\ldots,x\omega+y\bar\omega,\frac{x\omega+y\bar\omega}{1-\psi} \rangle_n.
$$
For example:
$$
\mathcal{L}^\omega (x,0; u,v) = J(x, 0; 0,u,v)
$$ 
For non-zero powers of $y$ the relations are slightly more complicated, and involve separating the $J$ function in a $\psi_\omega$ and a $\psi_{\bar\omega}$ part and applying variable shifts and integration to match the combinatorial factors. This so far has prevented us from finding a meaningful reformulation of our recursions in terms of the $J$ function. Of course, it would be very interesting if such a goal could be achieved.

\end{Rem}

\section{Table of Invariants of $\orb$}
\label{toi}
In the following table of invariants, we have set all the torus variables $t_i$ equal to $1$. We have boldfaced the non-equivariant invariants and decided to put a longer list of those.
\footnotesize
\begin{center}\renewcommand{\arraystretch}{1.25}
\begin{tabular}{|c|c|c|c|c|c|c|}\hline
$\rule{0ex}{3ex} n_1+n_2\downarrow$ & \multicolumn{5}{|c|}{$\langle\omega^{n_1}\overline{\omega}^{n_2}\rangle$ with $n_1\equiv n_2 \pmod 3$} \\ \hline
$\rule{0ex}{3ex} 3$ & $\mathbf{\frac 13}$ & $\frac 13$ &&& \\ \hline
$\rule{0ex}{3ex} 4$ & $-\frac 13$ & & & & \\ \hline
$\rule{0ex}{3ex} 5$ & $\frac 19$ & $\frac 29$ & & & \\ \hline
$\rule{0ex}{3ex} 6$ & $\mathbf{-\frac{1}{27}}$ & $-\frac{8}{27}$ & $-\frac{10}{27}$ & & \\ \hline
$\rule{0ex}{3ex} 7$ & $\frac{7}{27}$ & $\frac{19}{27}$ &&& \\ \hline
$\rule{0ex}{3ex} 8$ & $-\frac{5}{27}$ & $-\frac{98}{81}$ & $-\frac{179}{81}$ && \\ \hline
$\rule{0ex}{3ex} 9$ & $\mathbf{\frac 19}$ & $\frac{398}{243}$ & $\frac{1274}{243}$ & $\frac{686}{81}$ & \\ \hline
$\rule{0ex}{3ex} 10$ & $-\frac{451}{243}$ & $-\frac{905}{81}$ & $-\frac{6172}{243}$ && \\ \hline
$\rule{0ex}{3ex} 11$ & $\frac{1319}{729}$ & $\frac{14734}{729}$ & $\frac{52189}{729}$ & $\frac{100762}{729}$ & \\ \hline
$\rule{0ex}{3ex} 12$ & $\mathbf{-\frac{1093}{729}}$ & $-\frac{7684}{243}$ & $-\frac{400010}{2187}$ & $-\frac{38884}{81}$ & $-\frac{612100}{729}$ \\  \hline
$\rule{0ex}{3ex} 15$ & $\mathbf{\frac{119401}{2187}}$ &  &  &  &  \\  \hline
$\rule{0ex}{3ex} 18$ & $\mathbf{-\frac{27428707}{6561}}$ &  &  &  &  \\  \hline
$\rule{0ex}{3ex} 21$ & $\mathbf{\frac{102777653467}{177147}}$ &  &  &  &  \\  \hline
$\rule{0ex}{3ex} 24$ & $\mathbf{-\frac{210755831694887}{1594323}}$ &  &  &  &  \\  \hline
$\rule{0ex}{3ex} \lfloor\frac{n_2}{3}\rfloor\to$ & $0$ & $1$ & $2$ & $3$ & $4$ \\ \hline
\end{tabular}\end{center}
\normalsize

\addcontentsline{toc}{section}{Bibliography}
\bibliographystyle{alpha}
\bibliography{biblio}

\end{document}